\documentclass{gtart_h}  

\def\ifplaintex{\expandafter\ifx\csname documentclass\endcsname\relax}

\ifplaintex 
\else
\fi

\expandafter\ifx\csname beginpicture\endcsname\relax
\expandafter\ifx\csname documentclass\endcsname\relax
\input pictex \else
\input prepictex \input pictex \input postpictex \fi\fi

\def\gt{{\mathsurround=0pt\it $\cal G\mskip-2mu$eometry \&\ 
$\cal T\!\!$opology}}        

\def\gtp{{\mathsurround=0pt\it $\cal G\mskip-2mu$eometry \&\ 
$\cal T\!\!$opology $\cal P\!$ublications}}  

\def\lognumber#1{\def\thelognumber{#1}}
\def\volumenumber#1{\def\thevolumenumber{#1}}
\def\papernumber#1{\def\thepapernumber{#1}}
\def\volumeyear#1{\def\thevolumeyear{#1}}

\def\pagenumbers#1#2{\def\startpage{#1}\def\finishpage{#2}}
\def\published#1{\def\publishdate{#1}}
\def\proposed#1{\def\theproposer{#1}}
\def\seconded#1{\def\theseconders{#1}}
\def\received#1{\def\receiveddate{#1}}
\def\revised#1{\def\reviseddate{#1}}
\def\accepted#1{\def\accepteddate{#1}}

\def\coverauthors#1{\def\thecoverauthors{#1}}
\def\asciiauthors#1{\def\theasciiauthors{#1}}

\long\def\asciiabstract#1{\long\def\theasciiabstract{#1}}

\let\\\par\let\thelognumber\relax
\let\thevolumenumber\relax\let\thepapernumber\relax
\let\thevolumeyear\relax\let\thesamplenumber\relax\let\startpage\relax
\let\finishpage\relax\let\publishdate\relax\let\receiveddate\relax
\let\reviseddate\relax\let\accepteddate\relax\let\theasciititle\relax
\let\theasciiauthors\relax
\let\theasciiabstract\relax
\let\theasciiemail\relax\let\theshortauthors\relax\let\theshorttitle\relax
\let\thecoverauthors\relax

\long\def\maketitlep{   

\count0=\startpage

\gt\hfill      
\beginpicture
\setcoordinatesystem units <0.33truein, 0.33truein> point at 2.2 0.9
\setplotsymbol ({$\cal G$})
\plotsymbolspacing=9truept
\circulararc 315 degrees from 0 1 center at 0 0
\setplotsymbol ({$\cal T$})
\circulararc 315 degrees from 1 -1 center at 1 0
\endpicture
\break
{\small\ifx\thesamplenumber\relax 
Volume \else Sample
\fi\thevolumenumber\ (\thevolumeyear)
\startpage--\finishpage\nl
Published: \publishdate}
\vglue 0.5truein plus 0.4fil minus 0.1truein

{\parskip=0pt\leftskip 0pt plus 1fil\def\\{\par\smallskip}{\ifplaintex\large
\else\Large\fi\bf\thetitle}\par\medskip}   

\vglue 0pt plus 0.1fil 

{\parskip=0pt\leftskip 0pt plus 1fil\def\\{\par}{\sc\theauthors}
\par\medskip}

\vglue 0pt plus 0.1fil 

{\small\parskip=0pt\let\newline\\
{\leftskip 0pt plus 1fil\def\\{\par}{\sl\theaddress}\par}
\expandafter\ifx\theemail\relax    
\relax\else\vglue 5pt plus 0.02fil minus 2pt\def\\{\stdspace{\rm 
and}\stdspace} 
\cl{Email:\stdspace\tt\theemail}\fi
\ifx\theurl\relax                  
\relax\else\vglue 5pt plus 0.02fil minus 2pt\def\\{\stdspace{\rm 
and}\stdspace}
\cl{URL:\stdspace\tt\theurl}\fi\par}

\vglue 7pt plus 0.3fil minus 3pt

{\bf Abstract}
\vglue 5pt plus 0.1fil minus 2pt

\theabstract

\vglue 7pt plus 0.3fil minus 3pt

{\bf AMS Classification numbers}\quad Primary:\quad \theprimaryclass

Secondary:\quad \thesecondaryclass

\vglue 5pt plus 0.3fil minus 2pt

{\bf Keywords:}\quad \thekeywords

\vglue 10pt plus 0.5fil minus 5pt

{\small  Proposed: \theproposer\hfill Received: \receiveddate\nl
Seconded: \theseconders\hfill 
\ifx\reviseddate\relax                         
Accepted: \accepteddate                        
\else
Revised: \reviseddate                          
\fi}
\eject
}       

\let\maketitlepage\maketitlep
\let\maketitle\maketitlepage

\font\phead=cmsl9 scaled 950
\font=cmsl9 scaled 1050
\font\pnum=cmbx10 scaled 913
\font\lnum=cmbx10 
\font\pfoot=cmsl9 scaled 950
\font=cmsl9 scaled 1050
\ifplaintex
\headline{\vbox to 0pt{\vskip -4.5mm\line{\small\phead\ifnum
\count0=\startpage ISSN 1364-0380 (on line)
1465-3060 (printed) \hfill {\pnum\folio}\else\ifodd\count0\def\\{ }%
\ifx\theshorttitle\relax\thetitle\else\theshorttitle\fi\hfill{\pnum\folio}
\else\def\\{ and }{\pnum\folio}\hfill\ifx\theshortauthors\relax\theauthors
\else\theshortauthors\fi\fi\fi}\vss}}
\footline{\vbox to 0pt{\vglue 0mm\line{\small\pfoot\ifnum\count0=\startpage
\copyright\ \gtp\hfill\else
\gt, Volume \thevolumenumber\ (\thevolumeyear)\hfill\fi}\vss
}}
\else
\makeatletter
\def\@oddhead{{\small\lhead\ifnum\count0=\startpage ISSN 1364-0380 (on line)
1465-3060 (printed) \hfill {\lnum\number\count0}\else\ifodd\count0
\def\\{ }\ifx\theshorttitle\relax \thetitle \else\theshorttitle\fi\hfill
{\lnum\number\count0}\else\def\\{ and }{\lnum\number\count0}
\hfill\ifx\theshortauthors\relax 
\theauthors\else\theshortauthors\fi\fi\fi}}\def\@evenhead{\@oddhead}
\def\@oddfoot{\small\lfoot\ifnum\count0=\startpage\copyright\ \gtp\hfill\else
\gt, Volume \thevolumenumber\ (\thevolumeyear)\hfill\fi}
\def\@evenfoot{\@oddfoot}
\makeatother
\fi

\newwrite\gtoutfile
\long\gdef\makeheadfile{  
{\def\\{, }\def\s{ }
\immediate\openout\gtoutfile head.xxx
\immediate\write\gtoutfile{Proxy-for: \ifx\theasciiauthors\relax
\theauthors\else\theasciiauthors\fi\s<\ifx\theasciiemail\relax\theemail\else\theasciiemail\fi>}
\immediate\write\gtoutfile{\noexpand\\}
\immediate\write\gtoutfile{Authors: \ifx\theasciiauthors\relax
\theauthors\else\theasciiauthors\fi}
{\def\\{ }\immediate\write\gtoutfile{Title: \ifx\theasciititle\relax
\thetitle\else\theasciititle\fi}}
\immediate\write\gtoutfile{Subj-class: GT or SG or MG etc}
\immediate\write\gtoutfile{MSC-class: \theprimaryclass\ifx\thesecondaryclass\relax\else, \thesecondaryclass\fi}
\immediate\write\gtoutfile{Journal-ref: Geom. Topol. \thevolumenumber
(\thevolumeyear) \startpage-\finishpage}
\immediate\write\gtoutfile{Comments: Published by Geometry and Topology at}
\immediate\write\gtoutfile{\s\s http://www.maths.warwick.ac.uk/gt/GTVol\thevolumenumber/paper\thepapernumber.abs.html}
\immediate\write\gtoutfile{\noexpand\\}
\immediate\write\gtoutfile{}
\ifx\theasciiabstract\relax
\immediate\write\gtoutfile{\theabstract}\else
\immediate\write\gtoutfile{\theasciiabstract}\fi
\immediate\write\gtoutfile{}
\immediate\write\gtoutfile{\noexpand\\}
\immediate\write\gtoutfile{}
\immediate\closeout\gtoutfile}}  

\def\maketitlepage{\maketitlep\makeheadfile}
\let\maketitle\maketitlepage

\lognumber{362}

\received{20 September 2003}

\volumenumber{8}\papernumber{32}\volumeyear{2004}
\pagenumbers{1189}{1226}
\revised{22 August 2004}
\published{10 September 2004}
\accepted{11 July 2004}
\proposed{Leonid Polterovich}
\seconded{Yasha Eliashberg, Robion Kirby}

\usepackage{amsmath,amssymb}

\renewcommand{\Tilde}{\widetilde}
\newcommand{\R}{\mathbb{R}}

\newcommand{\Z}{\mathbb{Z}}
\newcommand{\C}{\mathbb{C}}

\newcommand{\CP}{\mathbb{CP}}
\newcommand{\p}{{\partial}}
\newcommand{\al}{{\alpha}}

\newcommand{\Om}{{\Omega}}
\newcommand{\om}{{\omega}}
\newcommand{\eps}{{\varepsilon}}
\newcommand{\de}{{\delta}}
\newcommand{\De}{{\Delta}}
\newcommand{\ga}{{\gamma}}
\newcommand{\Ga}{{\Gamma}}
\newcommand{\ka}{{\kappa}}
\newcommand{\la}{{\lambda}}
\newcommand{\si}{{\sigma}}
\newcommand{\graph}{{\rm graph}}
\newcommand{\Symp}{{\rm Symp}}
\newcommand{\Ham}{{\rm Ham}}
\newcommand{\Diff}{{\rm Diff}}

\newcommand{\Aa}{{\cal A}}

\newcommand{\Dd}{{\cal D}}

\newcommand{\Ll}{{\cal L}}
\newcommand{\Jj}{{\cal J}}

\newcommand{\Hh}{{\cal H}}

\newcommand{\Mm}{{\cal M}}

\newcommand{\Q}{\mathbb{Q}}

\newcommand{\La}{{\Lambda}}

\newcommand{\Si}{{\Sigma}}

\newcommand{\area}{{\rm area}}

\newcommand{\Totvar}{{\rm Totvar}}
\newcommand{\vol}{{\rm vol}}

\newcommand{\codim}{{\rm codim}}
\newcommand{\Crit}{{\rm Crit}}

\newcommand{\id}{{\rm id }}

\newcommand{\IMP}{{\Longrightarrow}}

\newtheorem{theorem}{Theorem}[section]
\newtheorem{cor}[theorem]{Corollary}

\newtheorem{prop}[theorem]{Proposition}

\newtheorem{lemma}[theorem]{Lemma}

\theoremstyle{definition}

\newtheorem{definition}[theorem]{Definition}

\newtheorem{conj}[theorem]{Conjecture}

\begin{document}

\title{A field theory for symplectic fibrations over surfaces}

\authors{Fran\c{c}ois Lalonde}
\coverauthors{Fran\noexpand\c{c}ois Lalonde}
\asciiauthors{Francois Lalonde}

\address{Department of Mathematics and Statistics, University of
  Montreal\\Montreal H3C 3J7, Quebec, Canada}

\email{lalonde@dms.umontreal.ca}

\begin{abstract}
We introduce in this paper a field theory on symplectic
manifolds that are fibered over a real surface with interior marked
points and cylindrical ends.  We assign to each such object a morphism
between certain tensor products of quantum and Floer homologies that are
canonically attached to the fibration. We prove a composition theorem
in the spirit of QFT, and show that this field theory applies naturally
to the problem of minimising geodesics in Hofer's geometry. This work
can be considered as a natural framework that incorporates both the
Piunikhin--Salamon--Schwarz morphisms and the Seidel isomorphism.
\end{abstract}

\asciiabstract{%
We introduce in this paper a field theory on symplectic manifolds that
are fibered over a real surface with interior marked points and
cylindrical ends.  We assign to each such object a morphism between
certain tensor products of quantum and Floer homologies that are
canonically attached to the fibration. We prove a composition theorem
in the spirit of QFT, and show that this field theory applies
naturally to the problem of minimising geodesics in Hofer's
geometry. This work can be considered as a natural framework that
incorporates both the Piunikhin-Salamon-Schwarz morphisms and the
Seidel isomorphism.}



\maketitle

\section{Introduction}

We describe in this paper a field theory for Hamiltonian fibrations
over oriented surfaces with boundaries and marked points. In this
context, a Hamiltonian fibration $(M,\om) \hookrightarrow P \to \Si $
is a symplectic manifold $(P,\Om)$ whose form restricts to a symplectic
form on each fiber, endowed with interior marked points on $\Si$. The
restriction of $\Om$ to the 1--skeleton of $\Si$ is assumed to be
fiberwise symplectically trivialisable, although such a trivialisation
is not part of the data. To each boundary component of $\Si$ there corresponds
a geometric version of Floer homology and  to each marked point the
quantum homology of the fiber at that point. Any partition in two subsets
of the boundary components and marked points gives rise to a morphism
between the corresponding tensor product of homologies. The gluing of
two such surfaces (using gluing on the boundary components and index zero
surgeries on a pair of marked points) and the compatible gluing of the
fibrations give rise to the composition of morphisms. This setting is a
geometric generalisation of the Piunikhin--Salamon--Schwartz \cite{PSS}
and the Seidel homomorphisms \cite{Se}.

We will show below that such a generalisation is natural in the study of
Hofer's geometry: we will use our field theory to give a new proof that a
symplectic isotopy $\phi_t \in \Diff_{\Ham}(M)$ generated by a Hamiltonian
$H_t$ is length minimising in the Hofer metric in its homotopy class
with fixed endpoints if two special classes corresponding to the fixed
minimum and maximum of the generating Hamiltonian are essential in the
Floer homology of $H_t$. The proof is geometric and works for all weakly
monotone manifolds (actually, it very likely holds for all symplectic
manifolds if one consider the virtual moduli space of pseudo-holomorphic
curves, but we will restrict here to the case where $(M,\om)$ is weakly
monotone). See Corollary~\ref{cor:Hofer} for new results on length
minimising geodesics in Hofer's geometry that this field theory yields.
Here is a simple example of application of Corollary~\ref{cor:Hofer}:
\medskip

{\bf Example}\qua Given a symplectic manifold with minimal Chern number
not in the range $[1,n]$, $\CP^n$ for instance, let $f\co M \to \R$ be
a Morse function.  Let $P$ be a global maximum of $f$ and $Q$ a global
minimum of $f$ such that the linearised flows have no nontrivial closed
orbit in time $\le 1$. Let $\eps > 0$ be sufficiently small so that
the graph of the diffeomorphism $\phi^f_t$ induced by $f$ be transversal
to the diagonal in $M \times M$ for all $t \in (0,\eps]$. Consider any
Hamiltonian $H_{t \in [0,1]}$ such that:
\begin{enumerate}
\item[(1)] on some open neighbourhoods of  $P$ and $Q$, $H_t = f$ for
  all $t \in [0,1]$
\item[(2)] $H_t = f$ for $t \le \eps$,
\item[(3)] the graph of $\phi^H_t$ remains transversal to the diagonal
  in $M \times M$ for all $t \in [\eps,1]$, and
\item[(4)] $H_t(x) \in [f(Q), f(P)]$ for all $x \in M$ and $t \in [0,1]$.
\end{enumerate}
Then the path induced by $H_{t \in [0,1]}$ is length minimising rel
endpoints in its homotopy class, in the Hofer metric. See the proof
of this claim after Corollary~\ref{cor:Hofer} in Section~\ref{se:EimpliesM}.

This paper is a natural extension of the ideas presented in
Lalonde--McDuff--Polterovich \cite{LMP} and in Entov \cite{E}. See the
remark at the end of Section~\ref{se:fieldtheory} for a generalisation
of this setting to a category of symplectic fibrations slightly larger
than the Hamiltonian ones.

We end this paper with a somewhat independent section treating a notion of
``total norm'' on Hamiltonian diffeomorphisms, defined as the sequence
of its genus--$g$ norms. The results and conjecture of that section turn
out to be closely related to Entov's paper  \cite{Entov} on the relation
of K--area and the commutator length of Hamiltonian diffeomorphisms.

\medskip
{\bf Acknowledgements}\qua I wish to thank the referee, whose questions
helped clarify some aspects of the first version.  This research was
partially supported by a CRC grant, NSERC grant OGP 0092913, and FCAR grant
ER--1199.

\section{A field theory for fibrations} \label{se:fieldtheory}

\subsection{Definition of the objects} \label{ss:objects}

Here is the definition of the basic objects. Let $\Si = \Si_{g,k,\ell}$ be
a real compact oriented surface of genus $g$ with interior marked points
$p_1,\ldots, p_k$ and oriented boundary components $S_1,\ldots,S_{\ell}$
(here $g$ is by definition the genus of the closed surface obtained by
capping with discs all boundary components). Let $(M, \om) \hookrightarrow
P \stackrel{\pi}{\to} \Si$ be a Hamiltonian fibration with fiber
{\em any} compact symplectic manifold (we will restrict ourselves
to simpler manifolds later when we will define $\Phi_{\si}$). By the
characterisation of such fibrations in \cite{MS, LM-Topology}, we may
assume that each fiber $M_b$ is equipped with a symplectic structure
$\om_b$ that admits an extension to a symplectic form $\Om$ on $P$
(note that the restriction of $\Om$ to each $M_b$ is equal to $\om_b$ as
forms, not only as cohomology classes). We may further assume that the
restriction of $P$ to the 1--skeleton  $\Si_1$ of $\Si$ (including its
boundary components) is fiberwise symplectically trivialisable. We
will denote by $\rho\co\pi^{-1}(\Si_1) \to \Si_1 \times M$ such a
trivialisation. Thus $\rho_{\ast} (\Om)$ restricts to $\om$ on each
fiber $\{b\} \times M$. Such a choice of trivialisation is not unique in
general at the homotopy level. The restriction of $\Om$ to $\pi^{-1}(S_j)$
has a one-dimensional kernel transversal to the fibers: the monodromy
of this foliation followed in the direction of the boundary component
corresponding to its orientation is a Hamiltonian diffeomorphism
$\phi_j\co M_{q_j} \to M_{q_j}$.

\begin{definition}
A {\em marked Hamiltonian symplectic fibration} $P$ over $\Si$ is given
by the following data:

\begin{enumerate}
\item[(1)] a fibration $(M,\om) \hookrightarrow (P,\Om) \to \Si$ over a
  real compact oriented surface with $\ell$ oriented boundary components,
  each one with a base point $q_j \in S_j$, together with $k$ interior
  marked points $p_1,\ldots, p_k \in \Si$;
\item[(2)] symplectic identifications $\eta_i\co (M_{p_i}, \om_{p_i}) \to
  (M, \om)$ for each $1 \le i \le k$ and $\eta_j\co (M_{q_j}, \om_{q_j})
  \to (M, \om)$ for each $1 \le j \le {\ell}$;
\item[(3)] a choice, for each $j$,  of an element $\xi_j$ of the
  universal cover $\Tilde{\Diff}_{\Ham}(M,\om)$ which projects to $\psi_j
  =_{def} \eta_j \circ \phi_j \circ \eta_j^{-1} \in \Diff_{\Ham}(M,\om)$;
\item[(4)] a partition in two subsets $\cal A_+, \cal A_-$ of the set
  $\{p_1,\ldots, p_k, S_1,\ldots,S_{\ell} \}$.
\end{enumerate}

Here the orientation ${\cal O}_j$ of the ``ingoing'' $S_j$s, ie
those that belong to $\cal A_+$, is such that ${\cal O}_j$ followed
by the inward normal orientation is equal to the orientation of the
surface. The orientation ${\cal O}_j$ of the ``outgoing'' $S_j$s, ie
those that belong to $\cal A_-$, is such that ${\cal O}_j$ followed by
the outward normal orientation is equal to the orientation of the surface.
\end{definition}

Note that, although the identifications $\eta$s are part of the data,
the choice of a trivialisation $\rho$ is not. In the sequel, any choice
of a trivialisation $\rho$ will  be compatible with the $\eta_j$s,
$1 \le j \le \ell$.

The main goal of this section is to assign to each end of $\Si$ a Floer
homology generated by the one-turn closed leaves of the characteristic
foliation of $\Om |_{\pi^{-1}(S_j)}$ and to each marked Hamiltonian
fibration $P \to \Si$ a morphism
\begin{multline*}
\Phi_{\si}\co \bigotimes_{p_i \in \cal A_+} QH_*(M_{p_i}) \otimes
  \bigotimes_{S_j \in \cal A_+} FH_*(S_j) \\
  \longrightarrow \bigotimes_{p_i \in \cal A_-} QH_*(M_{p_i}) \otimes
  \bigotimes_{S_j \in \cal A_-} FH_*(S_j)
\end{multline*}
depending on the choice of a homology class of sections $\si$, whose
domain is the tensor product over the ingoing data of the Floer homologies
(assigned to the ends) and the quantum homologies (assigned to the fibers
over the interior marked points) and whose codomain is the corresponding
tensor product over the outgoing data. This morphism is a generalisation
to fibrations of the morphisms introduced in Piunikhin--Salamon--Schwarz: it
is obtained by counting sections of $P$ that obey a certain Cauchy--Riemann
equation. However, in the counting of pseudo-holomorphic curves, we do
not fix the conformal structure of the base.  Note  however that all the
interesting cases that we have in mind involve marked surfaces $\Si$ with
cylindrical ends that admit, up to the data preserving diffeomorphism
group, only one conformal structure; therefore, the variation of the
conformal structure is not an issue in the applications that we have
in mind.

Here are some special cases of this construction.

\medskip
{\bf Example 1}\qua
When the base is $S^2$ with two marked points, the morphism:
$$\Phi_{\si}\co QH_*(M_{p_1}) \to QH_*(M_{p_2})$$
is the Seidel homomorphism (see \cite{Se,LMP}) that depends on the
homotopy class of the loop of Hamiltonian diffeomorphisms used to define
the fibration $P \to S^2$.

\medskip
{\bf Example 2}\qua
When the base is the $2$--disc with one interior marked point and the
fibration is trivial away from the end,  this is as in \cite{PSS} the usual
isomorphism between quantum and Floer homologies.

\medskip
{\bf Example 3}\qua
When the base is the $2$--sphere with three marked points and $P$ is the
trivial fibration (or three ends and the fibration is trivial away from
the cylindrical ends), this is the quantum product (or Floer product).

\medskip
It is not difficult to see that, when the base is the $2$--sphere with an
arbitrary number of ends and marked points in the ingoing data set and
only one element (either a point or a boundary component) in the outgoing
data set, our morphism is always a composition of the morphisms of the
last three examples, in the sense of our Composition Theorem (see below).
However, our theory makes sense in a wider context, where one may allow
monodromies over interior loops of $\Si$ to be non symplectically isotopic
to the identity -- see the last remark of this section.
   
The main result of this section is a Composition formula similar to
Theorem~3.7 of \cite{PSS}. It is useful even in the simplest case of
the gluing of two fibrations $(M, \om) \hookrightarrow P_i \to D^2$
(for $i=1,2$) over a disc with one interior marked point. We assume here that
the fibrations $P_1, P_2$  have the same Hamiltonian monodromy round
their end, and can therefore be glued together along the monodromies
to yield a new (not necessarily trivial) fibration $P_1 \# P_2$. The
Composition Theorem then states that the Seidel homomorphism
$$\Phi_{\si = \si_1 \# \si_2, P = P_1 \# P_2}\co QH_*(M) \to QH_*(M)$$
corresponding to $P_1 \# P_2$ factorises through the morphisms
$\Phi_{\si_1,P_1}\co QH(M) \to FH(\p P_1)$ and
$\Phi_{\si_2,P_2}\co FH(\p P_2) \to  QH(M)$. Thus, if one knows that
the Seidel homomorphism is an isomorphism, one concludes that both
$\Phi_{\si_1,P_1}$ and $\Phi_{\si_2,P_2}$ are isomorphisms too, and
conversely.  This will be used in our application to Hofer's geometry
in Section~\ref{se:EimpliesM}.

Let us return to the construction of the main objects. Our first goal
is to define a homology whose generators are the closed leaves of
the characteristic foliation of the restriction of $\Om$ to a given
boundary component, together with filling discs to be defined below,
that will be canonically identified with the Floer homology in a way
that is independent of the chosen trivialisation $\rho$.

Let $H_2^S(M,\Z)$ be the spherical part of the second homology group
of $M$, ie the image of the Hurewicz homomorphism $\pi_2(M) \to H_2(M,
\Z)$, and let $\La$ be the usual rational Novikov ring of the group $\Hh =
H_2^S(M,\Z)/\!\!\sim$  where $B\sim B'$
if $\om(B-B') = c(B-B') = 0$. Thus the elements of $\La$ have the form
$$\sum_{B\in \Hh} \la_B e^B$$
where for each $\ka$ there are only finitely many nonzero
$\la_B\in \Q$ with $\om(B) < \ka$.

Let us first quickly recall that the Floer homology $HF_*(H_t)$  can be
considered as the Novikov homology of the classical action functional
$$\Aa_{H}\co \Tilde \Ll \to \R$$
defined by $\Aa_H(x(t), v_x) = \int_{v_x} \om - \int_0^1 H_t(x(t)) dt$.
Here $\Tilde \Ll$ is a covering of the loop space $\Ll$ of all contractible
loops in $M$ defined in the following way.
Fix a constant loop $x_* \in \Ll$. Then the covering $\Tilde \Ll$ is defined by
requiring that two paths from $x_*$ to $x(t)$ are equivalent if the $2$--sphere
$S$ obtained by gluing the corresponding discs has $\om(S) = c(S) = 0$.
Thus the covering group
of $\Tilde \Ll$ is $\Hh$, and the Floer homology is a $\La$--module in
a natural way.

Now let us define the geometric Floer homology attached to an end
$S_j$. Note that any fiberwise symplectic trivialisation  $\rho_j\co
\pi^{-1}(S_j) \to S_j \times M$ pushes the form $\Om$ to a form $\Om'$
on $S_j \times M$ whose characteristic flow round the boundary projects on
the fiber $(M, \om)$ to a Hamiltonian isotopy $\psi_{t \in [0,1]}$. Pick
the trivialisation $\rho_j$ so that this isotopy belongs to the class
$\xi_j$. This is always possible by composing the trivialisation, if
needed, with an appropriate loop of Hamiltonian diffeomorphisms. Let
$P_j$ be the canonical Hamiltonian fibration over the 2--disc whose
monodromy round the boundary lies in the class $\xi_j$: this is given as
the part $P^-(H)$ under the graph of the Hamiltonian $H_{t \in [0,1]}$
that generates $\psi_{t \in [0,1]}$ -- see the precise definition in
Section~\ref{sse:hofer}. Glue $P_j$ to $P$ via the map $\rho_j$, and
denote by $\bar P$ the space obtained from $P$ by gluing to $P$ all the
$P_j$s, $1 \le j \le \ell$. We will call $\bar P$ the {\em closure of}
$P$. It is endowed with a natural ruled symplectic form $\bar{\Om}$ that
extends $\Om$; it extends as well the symplectic forms on the fibers of
$\bar P$. We will denote by the same symbol $P_j$ its image in $\bar P$.

\begin{definition}
\label{def:complex}
The {\em Floer complex} corresponding to $S_j$ has generators $(x(t),\!v)$
where $x$ is a parametrisation of a closed leaf of the characteristic
foliation (induced by a fixed parametrisation of the base $S_j$),
and where $v$ is a 2--disc in $P_j \subset \bar P$ with boundary
equal to $x$. Via the trivialisation $\rho_j$, these objects are in
one-to-one correspondence with the generators of the usual Floer complex
of $H_t$. Define the boundary operator as the pull-back through $\rho_j$
of the usual Floer boundary operator.
\end{definition}

\begin{lemma}
The Floer complex corresponding to $S_j$ depends only on the class $\xi_j$.
\end{lemma}

\begin{proof}
If $\rho_j'$ is another trivialisation in class $\xi_j$, the two
corresponding Hamiltonian isotopies $\psi_j$ and $\psi'_j$ are homotopic
with fixed endpoints. There is however a delicate point here: the
identification of these two isotopies (and therefore the identification
of the corresponding cylinders  $P_j, P_j'$, capping discs $v,v'$ and
boundary operators) depends on a choice of a homotopy between the two
paths $\psi_j$ and $\psi'_j$. But two different homotopies differ by an
element in $\pi_2(\Diff_{\Ham}(M))$ and it is shown in \cite{LM-Topology},
Proposition~5.4, that such elements act trivially on the homology of
$M$. Thus all these identifications act in the same way on capping
discs. It is then easy to see that they preserve the boundary operator
too. Note that this is equivalent to saying that the Seidel morphism
induced by a contractible loop $\phi_t$ of Hamiltonian diffeomorphisms
admits a unique lift to the universal cover of the loop space (independent
of the choice of the homotopy between the loop $\phi_t$ and the trivial
loop).
\end{proof}

We will denote by $FC(S_j)$ the Floer complex and by $FH(S_j)$ its
homology (or by $FC(S_j, \xi_j), FH(S_j, \xi_j)$ if we wish to emphasise
the dependence on $\xi_j$). It is a $\La$--module, and there is a canonical
$\La$--module isomorphism between $FH(S_j, \xi_j)$ and the Floer homology
of $M$.
Recall that this complex is graded by the Conley--Zehnder index normalised
so that it corresponds to the Morse index in the case of a $C^1$--small
autonomous Hamiltonian.

One may define the quantum complex associated with a Smale--generic pair
$(\kappa,g_M)$ consisting of a Morse function $\kappa\co M \to \R$ and a
metric $g_M$ on $M$: this is the usual Morse complex of $(\kappa,g_M)$
tensored with $\La$, ie $C_*(\kappa, g_M) \otimes \La$. When there
is no confusion, we will often omit the data $\kappa, g_M$ and  denote
$C_*(\kappa, g_M)$ by $C_*(M)$. The quantum complex is defined by $QC_*(M)
= C_*(M) \otimes \La$. The quantum homology of $M$ is defined as
$QH_*(M) = H_*(M)\otimes\La$,
$\Z$--graded with $\deg(a\otimes e^B) = \deg(a) - 2c(B)$. Associate to
each marked point $p_i$ the quantum homology of the fiber at that point,
$QH(M_{p_i})$.
Each $\eta_i$ identifies the quantum homology of  the fiber at $p_i$ with the quantum homology of $M$. 

\begin{lemma}
\label{le:sections}
Every Hamiltonian fibration over a closed oriented real surface admits
sections. Actually, the homology classes of sections form an affine
space over $H^S_2(M)$.
\end{lemma}

\begin{proof}
This has been proved in \cite{LMP} for Hamiltonian fibrations over the
$2$--sphere. Let $P \to \Si$ be a Hamiltonian fibration over an arbitrary
Riemann surface with boundary, given with a symplectic fiberwise
trivialisation over the $1$--skeleton $\Si_1$ of $\Si$. Thus finding a
section of $P$ which is constant over $\Si_1$ is equivalent to finding
sections of the quotient bundle over $\Si/\Si_1$. But the latter is a
bouquet of $2$--spheres, and the result therefore follows from the case
proved in \cite{LMP}.
\end{proof}

\subsection{Definition of the morphisms}
\label{ss:morphisms}
To simplify the exposition and avoid dealing with virtual moduli spaces,
assume from now on that $(M, \om)$ is weakly monotone. Recall that this
means that for every spherical homology class $B\in H_2(M)$
$$\om(B) > 0,\;\;c_1(B) \ge 3-n \quad\IMP \quad c_1(B)\ge 0.$$
This condition is satisfied if
either $\dim M \le 6$ or $M$ is {\em semi-monotone}, ie there  is
a constant $\mu \ge 0$ such that, for all spherical homology classes
$B\in H_2(M)$,
$$c_1(B) = \mu \om(B).$$
Weak monotonicity is the condition under which the ordinary (non virtual)
theory of $J$--holomorphic curves behaves well.

Define an equivalence relation on the set of homology classes of
sections of $\bar P$ by identifying two such classes if their values
under $c_{vert}$ and $\bar{\Om}$ are equal. Here $c_{vert}$ is the first
Chern class of the bundle $V \to \bar P$ whose fiber over $p \in \bar
P$ is the tangent space $T_{p}(\pi^{-1} (\pi(p)))$ to the fiber of the
$M$--bundle passing through $p$. Note that this definition is independent
of the extension $\bar{\Om}$ of $\om$.  Given an equivalence class of
sections $\si$ of $\bar P$,
whose existence is guaranteed by Lemma \ref{le:sections}, we are going
to define a morphism of chain complexes
$$\Phi_{\si}\co \bigotimes_{p_i \in \cal A_+} QC_*(M_{p_i}) \otimes
  \bigotimes_{S_j \in \cal A_+} FC_*(S_j)
\longrightarrow \bigotimes_{p_i \in \cal A_-} QC_*(M_{p_i}) \otimes
  \bigotimes_{S_j \in \cal A_-} FC_*(S_j).$$
that depends on the auxiliary data given by the choices of the generic
Smale--Morse pairs $(\kappa_i,g_i)$ on the fibers $M_{p_i}$, but will
turn out to be independent of these auxiliary data at the homology level.
Note that the domain and the codomain are both $\La$--modules in an obvious
way. Here the total degree of the domain or codomain is the sum of the
degrees of each factor, and we denote by $d$ the degree of $\Phi_{\si}$,
ie the integer such that $\deg (\Phi_{\si}(w)) = \deg (w) + d$ (see the
computation of $d$ below in Lemma \ref{le:index}). To simplify notations,
assume that the set $\cal A_+$ contains the indices from $i=1$ to $i=i_+$
and from $j=1$ to $j=j_+$. Any chain in the domain of $\Phi_{\si}$
decomposes as a sum of elements of the form: $(a_1 e^{\al_1} \otimes
\ldots \otimes a_{i_+} e^{\al_{i_+}})   \otimes ((x_1, v_1) \otimes \ldots
\otimes (x_{j_+}, v_{j_+}))$ where $a_1, \ldots, a_{i_+}$  are  unstable
manifolds of critical points of $\kappa_i$. The image of such an element
by $\Phi_{\si}$ has the form
\begin{multline*}
  \sum \; m \left(a_1, \ldots, a_{i_+}, b_{i_++1}, \ldots,
  b_k, (x_1,v_1), \ldots, (x_{\ell}, v_{\ell}), \si, \al_i\right) \\
  \left(b_{i_++1} e^{\al_{i_++1}} \otimes  \ldots \otimes
  b_k  e^{\al_{k}}\right) \otimes \left((x_{j_++1}, v_{j_++1}) \otimes \ldots
  \otimes (x_{\ell}, v_{\ell})\right)
\end{multline*}
where the summation runs through
all generators $(x_{j_++1}, v_{j_++1}), \ldots , (x_{\ell}, v_{\ell})$,
all sequences of 2--spheres $\al_{i_++1}, \ldots, \al_k$ and sequences of
stable manifolds $b_i$s associated with the $\kappa_i$s.  The number $m$
is defined to be equal to the cardinality, counted with signs, of the
moduli space
$$\Mm_{g,k,\ell}(a_1, \ldots, a_{i_+}, b_{i_++1}, \ldots,
  b_k, (x_1,v_1), \ldots, (x_{\ell}, v_{\ell}), \tau)$$
where
$$\tau = \si - \sum_{1 \le i \le i_+}\al_i
  + \sum_{i_+ + 1 \le i \le k}\al_i.$$
which is defined, roughly, by counting sections $u$ of the bundle that
are holomorphic away from the ends, satisfy a Floer-type equation on each
end, converge on each end to the loops $x_1, \ldots , x_{\ell}$, and
meet the $a$ and $b$ chains over the marked points of the surface. For
this, one needs in general to deal with a moduli space of conformal
structures on the base; however,  one cannot quotient out the pairs $(j,
u)$ made of a conformal structure and of a section by the diffeomorphism
group of $\Si$ that preserves the ends and the marked points, since the
action on the second factor does not yield a section. For this reason,
we will introduce the following definition:

\begin{definition}
Let $\Si$ be a real oriented surface of genus $g$ with $\ell$ ends
and $k$ interior marked points.  Let us denote by $\Jj_{g, \ell}$
the space of conformal structures on the interior of $\Si$ that are
standard on each end (ie biholomorphic to $\R/\Z~\times~[0, \infty)$)
and by $\Dd_{g,k, \ell}$ the group of compactly supported orientation
preserving diffeomorphisms of $\Si$ that fix pointwise the ends and
the marked points. Note that the quotient $\Jj_{g, \ell} / \Dd_{g,k,
\ell}$  has a natural smooth (orbifold) structure. We say that $\Si$ it
is {\em normalisable} if  the space of smooth sections of the quotient
map $\Jj_{g, \ell} \to \Jj_{g, \ell} / \Dd_{g,k, \ell}$  is non-empty
and connected. In such a case, we say that it is {\em normalised} if
the choice of a section has been made.
\end{definition}

Obviously, when the quotient space $\Jj_{g, \ell} / \Dd_{g,k, \ell}$
reduces to a point, which is the case when $\Si$ has genus $0$ and
the total number of ends and marked points is less than or equal to $3$,
the marked surface is normalisable: there is then only one class of
conformal structures to choose. In genus zero, with four marked points,
this quotient is the cross-ratio of the four points, and a normalisation
consists in  sending say $0,1, \infty$ to the first three points and
keeping the fourth one free.  In genus $g$ with zero marked points, a
normalisation is the choice of a particular structure for each class of
biholomorphic ones (ie a section of the Teichmuller space quotiented out
by the appropriate mapping class group). We will now assume that $\Si$
is normalised and will denote by $\si_{norm}$ the chosen section. It
is required that any complex structure in the image of $\si_{norm}$ be
standard on the cylindrical ends where the bump functions are defined;
it is allowed to vary everywhere else.

Here is a detailed description of this moduli space $\Mm_{g,k,\ell}$. In
$P$, there is a inner collar neighbourhood $V_j$ of the
hypersurface $W_j= \pi^{-1}(S_j)$  of the form $\pi^{-1}(U_j )$ where
$U_j=S^1 \times [0, \de]$. On this inner collar neighbourhood,
$\Om$ has the form $\Om_{W_j} + d(y d \theta) = \Om_{W_j} + dy \wedge d \theta$
where $y$ is the coordinate in $[0, \de]$. (Here $\de$ corresponds to
the boundary $S_j$ and $0$ is in the interior of the surface.) Therefore,
the characteristic
foliation of the hypersurfaces $W_j(y)$ does not depend on $y$: $W_j$ is {\em
stable} and actually flat. Take a generic family of almost
complex structures $J_b$ on each fiber $M_b$ of $P$, compatible with
$\om_b$, and let $\Hh$ be the symplectic connection which at $p \in P$ is
the $\Om$--orthogonal complement of the tangent space to the fiber at $p$.  
Consider the function
$h_j\co [0, \de] \to \R^+$  that vanishes near $0$, is non-decreasing
and reaches the slope $1$ at $ y = \de$, and let us denote by $f_j$ its
pull-back by the projection $S^1 \times [0,\de] \to [0,
\de]$. Thus $f_j$ has a periodic orbit of period $1$ at $\de$ with respect
to the symplectic form $dy \wedge d \theta$. Endowing  the $j$th end
$S^1 \times [0, \de]$ with the form $dy \wedge d \theta$ and the standard
conformal structure $ J_{st} \p_y = \p_{\theta}$, solves the equation
$$\p_{s} v + J_{st} \p_{t} v = - \nabla f_j$$
for a map $v\co S^1 \times [0, \infty) \to S^1 \times [0, \de)$ with
initial condition $v(t,0) = (t,0)$ which converges to the periodic orbit
of $f_j$ at $\ y = \de$. Push forward the conformal structure of the
domain of $v$ to its codomain: this gives a conformal structure $J_j$
on the $j$th end, conformally equivalent to the standard structure on
the semi-infinite cylinder. Now  each $J_{\Si} \in  \si_{norm}(\Jj_{g,
\ell} / \Dd_{g,k, \ell})$ is required to coincide with $J_j$ on each
end. For any such structure,  consider the
sections with finite energy $u\co\Si \to P$  which  are solutions of
the equation:
\begin{equation}
\renewcommand{\theequation}{$*$}
AC_{J_{\Si}, J_b} (\pi_{\Hh} \circ du) (b) = 0 
\end{equation}
for all $b \in \Si$, 
where $ \pi_{\Hh}$ is the projection on the tangent space to the fiber
defined by the connection $\Hh$ and $AC$ is the anti-complex part of
the derivative. Impose the following boundary conditions:
\begin{enumerate}
\item[(1)] $u(p_i)$ meets the chain $\iota(a_i)$ for $i \le i_+$ and
  $\iota(b_i)$ for $i > i_+$ (where $\iota$ denotes the map induced by
  inclusion of the fiber to $P$);
\item[(2)] $u(s,t) \to x_j(t)$ as $s \to \de$ on each end, and
\item[(3)] the section class that $u$ represents in $H_2(\bar P; \Z)$,
  once capped in the obvious way at each end by the $v_j$s $1 \le j \le
  \ell$, is equal to  $\tau$.
\end{enumerate}
\setcounter{equation}{0}
\renewcommand{\theequation}{\arabic{equation}}

Denote by 
$$\Mm^{\Ga}_{g,k,\ell,\si_{norm},\{J_b\},F } (a_1, \ldots, a_{i_+},
  b_{i_++1}, \ldots, b_k, (x_1,v_1), \ldots, (x_{\ell}, v_{\ell}), \tau)$$ 
the space of these solutions ($\Ga$ stands for sections), in which
$J_{\Si}$ is allowed to vary in the image of the chosen section of
$\Jj_{g, \ell} \to \Jj_{g, \ell} / \Dd_{g,k, \ell}$.   Here $F$ is the sum
of the $F_j\co P \to \R$ defined to be zero away from the ends and given,
on the above collar neighbourhood $V_j$, by $f_j  \circ \pi$ where $\pi$
is the projection of the fibration.

Thus this moduli space is the set of all pairs $(u, J_{\Si})$ where
$J_{\Si}$ belongs to the image of the generic section $\si_{norm}$
and $u$ is a section of the fibration $P \to \Si$ satisfying $(*)$ and
the above three boundary conditions. Because the total space of the
manifold is assumed to be weakly monotone, the  standard transversality
and orientation results  can be used to show that such a moduli space
is an oriented manifold for generic data  $\si_{norm},\{J_b\},F$ --
see for instance McDuff--Salamon \cite{MS2} and Salamon \cite{S} for
closed curves and Floer flow lines, and Gatien--Lalonde \cite{GL} for
a rigorous treatment of the transversality issues when the conformal
structure of the source $\Si$ is allowed to vary in a simple way.

{\bf Remark}\qua
Observe that a family of almost-complex structures on the fibers
$\{J_b\}$ and a complex structure $J_{\Si}$ on $\Si$ uniquely determine
an almost-complex structure  $J$ on $P$, compatible with $\Om$,  that
restricts to $\{J_b\}$ on the fibers,  is such that the projection
$\pi$ is $(J_{\Si},J)$--holomorphic, and has the horizontal distribution
$\Hh$ (ie the distribution of planes $\Om$--orthogonal to the fibers)
as invariant subspaces. It is then obvious that each above solution
corresponds to a section $u$ which is a solution of the equation
\begin{equation}
\renewcommand{\theequation}{$**$}
\bar{\p}_{J_{\Si}, J} u = -  \nabla F 
\end{equation}
where    $\nabla$ is the Riemannian gradient corresponding to the metric  $\Om(\cdot, J \cdot)$ and $F$ is defined as above. 
Note that we will not have to consider more general $J$s than the ones corresponding to a pair $(J_{\Si},\{J_b\})$ in this way.
\setcounter{equation}{0}
\renewcommand{\theequation}{\arabic{equation}}

Now,  define the number $m (a_1, \ldots, a_{i_+}, b_{i_++1}, \ldots,
b_k, (x_1,v_1), \ldots, (x_{\ell}, v_{\ell}), \tau)$ as zero if the
above moduli space has dimension different from $0$, and by the number
of elements of that moduli space (counted with sign) if its dimension
is zero.

\begin{prop}
\label{prop:chain}
The homomorphism $\Phi_{\si}$ is a morphism of chain complexes.
\end{prop}

\begin{proof}
Since we work in a weakly monotone manifold, this proposition can
be proved by standard techniques. Note first that the space $P$ is a
compact symplectic manifold with flat boundary. Thus the Gromov--Floer
compactness theorem applies as in the case of the standard Floer theory.

Each element of the moduli space
$$\Mm^{\Ga}_{g,k,\ell,\si_{norm},\{J_b\},F }  (a_1, \ldots, a_{i_+},
b_{i_++1}, \ldots, b_k, (x_1,v_1), \ldots, (x_{\ell}, v_{\ell}), \tau)$$
consists of three pieces:
\begin{enumerate}
\item[(1)] the ingoing flow lines from the $a_i$s ($i \le i_+$) to $u$,
\item[(2)] $u$ itself, and
\item[(3)] the outgoing flow lines from $u$ to the $b_i$s ($i > i_+$).
\end{enumerate}
Now consider the same moduli space where
exactly one of the terms in the tensor product of the elements in, say,
the target complex is replaced by a term appearing in its differential
(this amounts, by the Leibniz rule, to replacing the data in the target
space $$b_{i_++1}, \ldots, b_k,(x_{j_++1},v_{j_++1}), \ldots, (x_{\ell},
v_{\ell})$$ by a non-vanishing element appearing in its differential).
The index is then increased by one, hence this moduli space has real
dimension $1$ and has boundaries that are due either to the boundary
of the Morse chain complex associated to the quantum complexes, or to
the boundary of the moduli space of holomorphic sections. The former
consists of
broken flow lines of the Morse--Smale complexes $(\kappa_i,g_i)$
  where the breaking point is either in a fiber corresponding to the
  ingoing data, or in a fiber corresponding to the outgoing data (ie they
  correspond to the boundary of the stable or unstable chains);
and the latter consists of
broken flow lines of the Floer equations over one of the
  ingoing or outgoing ends.
By standard gluing techniques, either of these two events can be realised
as the boundary of a moduli space of the above form. Thus each element of
$$\Phi_{\si} \p \left( (a_1 e^{\al_1} \otimes \ldots \otimes a_{i_+}
  e^{\al_{i_+}})   \otimes ((x_1, v_1) \otimes \ldots \otimes (x_{j_+},
  v_{j_+})) \right)$$ 
can be paired, after cancellations,  with an element in  
$$\p \Phi_{\si} \left( (a_1 e^{\al_1} \otimes \ldots \otimes a_{i_+}
  e^{\al_{i_+}})   \otimes ((x_1, v_1) \otimes \ldots \otimes (x_{j_+},
  v_{j_+})) \right).$$ 
By the compactness theorem, there are for each
given energy level $\ka$  only finitely many  homology classes $D$ in $P$
with $\Om(D-\si) \le \ka$ that are represented by $J$--holomorphic curves
in $P$ with the fixed boundary conditions.  Thus
$\Phi_{\si}$ satisfies the finiteness condition for
elements  of the codomain. Finally, because the sum of a section
with a sphere in the fiber gives a homology class which is obviously
independent of the choice of the fiber where the sphere is taken (a fiber
in the domain of $\Phi_{\si}$ or in its codomain), the map $\Phi_{\si}$
is $\La$--linear.
\end{proof}

Thus the morphism $\Phi_{\si}$ descends to a homomorphism at the homology level. We omit the proof of the following proposition that follows by standard cobordism arguments as in Floer's theorem of invariance under change of auxiliary data:

\begin{prop} The map 
$$\Phi_{\si}\co \bigotimes_{p_i \in \cal A_+} QH_*(M_{p_i}) \otimes
  \bigotimes_{S_j \in \cal A_+} FH_*(S_j)
  \longrightarrow \bigotimes_{p_i \in \cal A_-} QH_*(M_{p_i})
  \otimes \bigotimes_{S_j \in \cal A_-} FH_*(S_j).$$
is well-defined: it is independent of all auxiliary choices that have
been made, ie it is independent of the choices of the fiberwise almost
complex structures $\{J_b\}$ on $M_b$, of the bump function $F$ and of
the (generic) choice of the normalisation.
\end{prop}

\begin{lemma}
\label{le:index}
The degree $d$ of $\Phi_{\si}$ is:
\begin{enumerate}
\item[\rm(0)] $d= 2 c_{vert}(\si)$ in the genus $0$ case;
\item[\rm(1)] $d = 2 c_{vert}(\si) -2n +2$ in the genus $1$ case;
\item[\rm(2)] $d = 2 c_{vert}(\si) + g(6-2n) - 6$ in the genus $g \ge 2$ case.
\end{enumerate}
\end{lemma}

\begin{proof}
The index formula for closed parametrised $J$--holomorphic curves $u$
of genus $g$ in an almost-complex  manifold $V$ of dimension $2n$ is:
$$\mbox{index} =  2(c_1(u) + n(1-g)) +  \dim {\cal T}_g$$
where the dimensions are real, ${\cal T}_g$ is the Teichmuller space and
$c_1$ denotes the first Chern class of the tangent bundle of $V$. Thus
the dimension of the space of $J$--holomorphic sections of $\bar P$ is
given by the same formula, where $c_1$ is replaced by $c_{vert}$. Hence
the dimension of the moduli space
$$\Mm^{\Ga}_{g,k,\ell,\si_{norm},\{J_b\},F} (a_1, \ldots, a_{i_+},
  b_{i_++1}, \ldots, b_k, (x_1,v_1), \ldots, (x_{\ell}, v_{\ell}), \tau)$$
is equal, for genus $0$, to 
\begin{multline*}
\sum_{1 \le j \le j_+} \mu(x_j,v_j) \; - \;  \sum_{j_+ + 1 \le j \le
  \ell} \mu(x_j,v_j) \; + \;  2 c_{vert}(\tau) \; + \;  2n \\
  - \sum_{1 \le i \le i_+} \codim \, a_i
  \; - \sum_{i_++1 \le i \le k} \codim \, b_i.
\end{multline*}
Thus the dimension vanishes when
\begin{multline*}
\sum_{i_++1 \le i \le k} \codim \, b_i =     \sum_{1 \le j \le j_+}
  \mu(x_j,v_j) \; - \; \sum_{j_+ + 1 \le j \le \ell} \mu(x_j,v_j) \\
  + \;  2 c_{vert}(\tau) \; + \; 2n - \sum_{1 \le i \le i_+} \codim \, a_i,
\end{multline*}
that is to say when 
\renewcommand{\theequation}{$*\!*\!*$}
\begin{multline}
\sum_{i_++1 \le i \le k} \dim a_i  = \sum_{1 \le j \le j_+} \mu(x_j,v_j)
  \; - \sum_{j_+ + 1 \le j \le \ell} \mu(x_j,v_j) \\
  + \;  2 c_{vert}(\tau)
  \; + \; 2n \; - \; \sum_{1 \le i \le i_+} \codim \, a_i.
\end{multline}
But \def\strut{\vrule width 0pt height 15pt depth 5pt}
$$\begin{array}{lll}
\deg\Phi_{\si}  =   \deg \Psi_{\si}( a_1 e^{\al_1} \otimes \ldots \otimes
  a_{i_+} e^{\al_{i_+}}  \otimes (x_1, v_1) \otimes \ldots \otimes (x_{j_+},
  v_{j_+})) && \\\strut
  - \deg (a_1 e^{\al_1} \otimes \ldots \otimes a_{i_+} e^{\al_{i_+}}   \otimes
  (x_1, v_1) \otimes \ldots \otimes (x_{j_+}, v_{j_+})) && \\\strut
= (\sum_{i_++1 \le i \le k} \dim a_i  - 2\sum_{i_++1 \le i \le k}
  c_{vert}(\al_i) + \sum_{j_++1 \le j \le \ell} \mu(x_j,v_j) ) && \\\strut
- ( \sum_{1 \le i \le i_+} \dim a_i  - 2 \sum_{1 \le i \le i_+}
  c_{vert}(\al_i) + \sum_{1 \le j \le j_+} \mu(x_j,v_j) ) && \\\strut
= ( \sum_{1 \le j \le j_+} \mu(x_j,v_j) - \sum_{j_+ + 1 \le j \le \ell}
  \mu(x_j,v_j) + 2 c_{vert}(\tau) + 2n && \\\strut
- \sum_{1 \le i \le i_+} \codim \, a_i - 2\sum_{i_++1 \le i \le k}
  c_{vert}(\al_i) + \sum_{j_++1 \le j \le \ell} \mu(x_j,v_j) ) && \\\strut
- ( \sum_{1 \le i \le i_+} \dim a_i - 2 \sum_{1 \le i \le i_+}
  c_{vert}(\al_i) + \sum_{1 \le j \le j_+} \mu(x_j,v_j) ) && \\\strut
= ( \sum_{1 \le j \le j_+} \mu(x_j,v_j) - \sum_{j_+ + 1 \le j \le \ell}
  \mu(x_j,v_j) + 2 c_{vert}(\si) + 2n && \\\strut
- \sum_{1 \le i \le i_+} \codim \, a_i + \sum_{j_++1 \le j \le \ell}
  \mu(x_j,v_j) ) - \sum_{1 \le i \le i_+} \dim a_i && \\\strut
- \sum_{1 \le j \le j_+} \mu(x_j,v_j)) && \\\strut
=   2 c_{vert}(\si) + 2n - \sum_{1 \le i \le i_+} \codim \, a_i -
\sum_{1 \le i \le i_+} \dim a_i = 2 c_{vert}(\si). &&
\end{array}$$
where the first equality relies on the definition of the degree of a
homomorphism, the second on the definition of the total degree of an
element, the third  on $(***)$, the fourth is obtained by replacing $\tau$
by its value $ \si - \sum_{1 \le i \le i_+} \al_i  + \sum_{i_+ + 1 \le
i \le k} \al_i$ and suppressing the terms that cancel out, and the last
two are straightforward computations. The computations in higher genera
are similar.
\end{proof}
\renewcommand{\theequation}{\arabic{equation}}
\setcounter{equation}{0}

Observe finally that
$$\Phi_{\si + B} = \Phi_{\si}\otimes e^{-B}.$$

\subsection{The composition theorem}
\label{ss:composition}

\begin{definition}
Let $P' \to \Si'=\Si'_{g',k',\ell'}$ and $P'' \to
\Si''=\Si''_{g'',k'',\ell''}$ be two marked Hamiltonian fibrations with
same fiber $(M, \om)$. Suppose that there is a bijection
$$\mu\co \{i'_+ +1, \ldots, k'\} \to \{1, \ldots, i''_+\}$$
between the ``outgoing'' $p'_i$s and the ``ingoing'' $p''_{i}$s and
a bijection
$$\nu\co \{j'_+ +1, \ldots, \ell'\} \to \{1, \ldots, j''_+\}$$
between the ``outgoing'' $S'_j$s and the ``ingoing'' $S''_{j}$s.
Assume moreover that the corresponding monodromies coincide, ie that
$\xi'_j = \xi''_{\nu(j)}$.

{\rm (1)}\qua The {\em gluing of $P'$ and $P''$}, denoted by $P=P' \# P''$,
  is by definition the marked Hamiltonian fibration obtained in the
  following way.  First extend the symplectic trivialisations $\eta'_i, \;
  i \in \{i'_+ +1, \ldots, k'\}$ over small discs where the form $\Om'$
  can therefore be identified with $(D^2 \times M,  \om_{st} \oplus
  \om)$ where $\om_{st}$ is the area form on the disc. Do the same near
  $\eta''_i, \;  i \in \{1, \ldots, i''_+\}$. Use these identifications
  to perform the obvious  surgery between $(\eta'_i)^{-1}(D^2 \times
  M)$ and $(\eta''_{\mu(i)})^{-1}(D^2 \times M)$ that covers the
  index $0$ surgery between $\Si'$ and $\Si''$ at the points $p'_i$
  and $p''_{\mu(i)}$. Similarly, consider a diffeomorphism $r_j\co S'_j \to
  S''_{\nu(j)}$ preserving both the orientations and the base points $q'_j,
  q''_{\nu(j)}$, and define
  $$\theta_j\co W'_j \to W''_{\nu(j)}$$
  as the symplectic diffeomorphism whose restriction to the fiber
  $M_{q'_j}$ is the map
  $$(\eta'')^{-1}_{\nu(j)}  \circ \eta'_j\co  M_{q'_j} \to M_{q''_{\nu(j)}}$$
  and which sends each fiber $M_b$ to its image $M_{r_j(b)}$ so that
  the characteristic foliations be preserved.

  Then glue the hypersurface $W'_j \to S'_j$ to the hypersurface
  $W''_{\nu(j)} \to S''_{\nu(j)}$ via the attaching map $\theta_j$. The
  resulting marked Hamiltonian fibration $P=P' \# P''$ over $\Si =
  \Si' \# \Si''$ has $\cal A'_+$ as ingoing data and $\cal A''_-$ as
  outgoing data.

{\rm (2)}\qua For each $i \in \{i'_+ +1, \ldots, \ell'\}$, there is an
  identification $\Phi_{QH}$ of $QC(M_{p'_i})$ with $QC(M_{p''_{\mu(i)}})$
  induced by the identification of the fibers.

  For each $j \in \{j'_+ +1, \ldots, \ell'\}$, there is an identification
  $\Phi_{FH}$ of $FC(S'_j)$ with $FC(S''_{\nu(j)})$ that maps $(x',v')$
  to an element $(x'',v'')$ in the following way. The closed leaf $x'(t)$
  is sent to $\theta(x'(t))$ and the disc $v' \subset P'_j$ is sent to
  the unique disc $v''$ in $P''_{\nu(j)}$ such that the gluing $v' \#
  v''$ in the fibration $P'_j \# P''_{\nu(j)} \to S^2$ is flat.
  To simplify notations, we will denote by $\Phi$ any tensor product of
  copies of $\Phi_{QH}$ and $\Phi_{FH}$.

{\rm (3)}\qua Finally, if $\si', \si''$ are sections of $\bar P', \bar{P''}$,
  we will denote by $\si' \#\si''$ the section of $\overline{P' \# P''}$
  such that in each $P'_j \# P''_{\nu(j)}$ the resulting sphere section
  be flat.

\end{definition}
  
Here is the main result of this section.

\begin{theorem}
\label{thm:composition}
Let $P'$ and $P''$ be two marked Hamiltonian fibrations and $P = P'
\#P''$ their gluing. Assume, to simplify the argument, that the number
of outgoing data of $P'$ is $1$.  Let $\si', \si''$ be section classes
in $P',P''$ and $\si = \si' \# \si''$ their gluing. Then
$$\Phi''_{\si''} \circ  \Phi \circ \Phi'_{\si'} = \Phi_{\si}$$
\end{theorem}

\begin{proof}[Sketch of the proof]\qua

(A)\qua To make the argument clearer, suppose first that each Morse chain
and periodic orbit at the source and the target spaces are {\em cycles}
and actually assume even that there is no element in their differential.
We will examine the general case of the construction of a chain homotopy
in (B) below.

First note that we may assume that the Morse--Smale data $(\kappa'_i,
g'_{M_i})$ at the point $p'_i$ coincides with the Morse--Smale data at
the point $p''_{\mu(i)}$.  So under these hypotheses, the statement
reduces to establishing the formula
\renewcommand{\theequation}{$*$}
\begin{multline}
\sum m_{g',k',\ell',\si'_{norm},\{J'_b\},F' }  \left(a'_1, \ldots,
  a'_{i'_+}, b'_{i'_++1}, \ldots, b'_{k'}, (x'_1,v'_1), \ldots, (x'_{\ell'},
  v'_{\ell'}), \tau'\right)  \\
\times m_{g'',k'',\ell'',\si''_{norm},\{J''_b\},F'' }
  \left(\Phi_{QH}(b'_{i'_++1}),
  \ldots, \Phi_{QH}(b'_{k'}), b''_{i''_+ + 1}, \ldots, b''_{k''},\right. \\
\hfill\left. \Phi_{FH}(x'_{j'_++1}, v'_{j'_++1}), \ldots,  \Phi_{FH}(x'_{\ell'},
  v'_{\ell'}),   (x''_{j''_+ +1},v''_{j''_+ +1}), \ldots,(x''_{\ell''},
  v''_{\ell''}), \tau''\right) \\
= m_{g,k,\ell,\si_{norm},\{J_b\},F} \left(a'_1, \ldots, a'_{i'_+},b''_{i''_+ +
  1}, \ldots, b''_{k''}, (x'_1,v'_1), \ldots, (x'_{j'_+}, v'_{j'_+}), \right. \\
\left.  (x''_{j''_++1}, v''_{j''_++1}), \ldots, (x''_{\ell''},
  v''_{\ell''}), \tau\right) 
\end{multline}
where the sum is taken over all elements $b'_{i'_++1}, \ldots, b'_{k'}$
of a basis of $H_*(M) \otimes \ldots \otimes H_*(M)$, all tensor
products $x'_{j'_++1} \otimes \ldots \otimes x'_{\ell'}$ of 1--turn closed
characteristic leaves, and all classes $\tau', \tau''$ such that $\tau'
\# \tau'' = \tau$.
\setcounter{equation}{0}
\renewcommand{\theequation}{\arabic{equation}}

The quantum part of the summation is taken over all elements
$b'_{i'_++1}, \ldots, b'_{k'}$ of a basis of $H_*(M) \otimes \ldots
\otimes H_*(M)$ corresponding to the outgoing data
$p'_i$ of the fibration $P'$ (and their images by $\Phi_{QH}$ in the
ingoing data of the second fibration $P''$). Since this is the basis
of the Morse complex corresponding to $(\kappa'_i, g'_{M_i})$ at the
$M$--fiber over the point $p'_i$, and because the boundary conditions of
the defining equation of the moduli space tell us that a solution $u'$
meets the {\em stable} submanifold of the critical points $b'_i$ while a
solution $u''$ meets the {\em unstable} submanifold of the corresponding
points, the above summation corresponds to taking the inverse image
${\cal I}$ by the evaluation map
\begin{multline*}
\Mm^{\Ga}_{g',k',\ell',\si'_{norm},\{J'_b\},F'} \left(a'_1, \ldots,
a'_{i'_+},  (x'_1,v'_1), \ldots, (x'_{\ell'}, v'_{\ell'}), \tau'\right)
\times \\
 \Mm^{\Ga}_{g'',k'',\ell'',\si''_{norm},\{J''_b\},F''} \left( b''_{i''_+ +
 1}, \ldots, b''_{k''}, \Phi_{FH}(x'_{j'_++1}, v'_{j'_++1}),  \ldots, \right. \\
\left. \Phi_{FH}(x'_{\ell'}, v'_{\ell'}),   (x''_{j''_+ +1},v''_{j''_+
 +1}), \ldots,(x''_{\ell''},
  v''_{\ell''}), \tau''\right) \\
 \longrightarrow \Pi_{i \in \{i'_++1, \ldots, , k' \}}  (M_i \times M_i)
\end{multline*}
of the product over the $i \in \{i'_++1, \ldots, k' \}$ of the cycles
$\sum_{r \in \Crit (\kappa_i)} b'_{i,r} \otimes a'_{i,r}$ where
$b'_{i,r}$ is the stable submanifold corresponding to a critical point
$q_r$ of the Morse function $\kappa'_i$ and $a'_{i,r}$ is the unstable
submanifold corresponding to the same critical point. But $\sum_{r \in
\Crit(\kappa_i)} b'_{i,r} \otimes a'_{i,r}$ is a cycle in $M_i \times M_i$
homologous to the diagonal  $\De_i \subset M_i \times M_i$. Hence the
proof reduces to establish the above equation $(*)$ when the intermediate
quantum data belong to the diagonal.

This boils down to a gluing theorem that gives a map from pairs of
solutions in ${\cal I}$ to solutions in the moduli space
\begin{multline*}
\Mm^{\Ga}_{g,k,\ell,\si_{norm},\{J_b\},F} \left(a'_1, \ldots,
  a'_{i'_+},b''_{i''_+ + 1}, \ldots, b''_{k''}, (x'_1,v'_1), \ldots,
  (x'_{j'_+}, v'_{j'_+}),\right. \\
\left. (x''_{j''_++1}, v''_{j''_++1}), \ldots, (x''_{\ell''},
  v''_{\ell''}), \tau \right).
\end{multline*}
The well-known gluing techniques in weakly monotone manifolds give such a
map: indeed, the gluing techniques for solutions near a fiber $M_{p_i}$
are exactly those exposed in Lalonde--McDuff--Polterovich \cite{LMP} and
developed in McDuff \cite{McD-sequel} (see also Ruan--Tian \cite{RT}),
while the gluing techniques for solutions on the asymptotic ends are
Floer's theorems on the composition of homotopies. The essential point
is to show that this gluing map is algebraically onto. This is what we
now prove.

Observe first that if $u$ is a solution in
\begin{multline*}
\Mm^{\Ga}_{g,k,\ell,\si_{norm},\{J_b\},F}  =
\Mm_{g,k,\ell,\si_{norm},\{J_b\},F} \left(a'_1, \ldots,
  a'_{i'_+},b''_{i''_+ + 1}, \ldots, b''_{k''}, \right. \\
\left. (x'_1,v'_1), \ldots,  (x'_{j'_+}, v'_{j'_+}),
  (x''_{j''_++1}, v''_{j''_++1}), \ldots, (x''_{\ell''},
  v''_{\ell''}), \tau \right)
\end{multline*}
then the $\Om$--area of $u$ is determined by the
boundary conditions and the homology class $\tau$, ie it is fixed
by the homology class of $u$ relative to its boundary conditions.

\begin{lemma}
\label{le:energy}
Let $u$ be a solution in $\Mm^{\Ga}_{g,k,\ell,\si_{norm},\{J_b\},F}
$. Then the $\Om$--area of each end is non-negative.
\end{lemma}

\begin{proof}
Recall that the restriction of a solution $u$ near the $j$th end is a map
$$u_j\co S^1 \times [0, \infty) \to P$$
that satisfies the equation $(**)$.
The computation of the
$s$--energy of the restriction of the solution  to that end  leads to
$$0 \le \int_{S^1 \times [0, \infty]}
  \left\| \frac{\partial u_j}{\partial s} \right\|^2 \le
  \left( \int_{S^1 \times [0, \infty]} (u_j)^{*}(\Om) \right)  - T_j$$
where $T_j > 0$ is the total variation of the function $F_j$. Thus the
$\Om$--area of that end is bounded from below by $T=T_j$, and is therefore
non-negative.
\end{proof}

Now let $u_0$ be a solution in $\Mm^{\Ga}_{g,k,\ell,\si_{norm},\{J_b\},F}$
and $\ga_j$ be a loop on $\Si$ in the class where the gluing of $S'_j$
with $S''_{\nu(j)}$ took place -- to simplify notation, we will now
drop the indices. Denote by $W$  the inverse image $\pi^{-1}(\ga)$;
note that the kernel of $\Om_W = \Om |_W$ is transverse to the
fibers. Recall that in a neighbourhood $U = \ga \times [-\de, \de]$ of
$\ga$, $\Om$ has the form $\Om'_W + d(s \al)$ where $\al$ is any $1$--form
on $W$ that does not vanish on the kernel of $\Om_W$ and $s$ is the
coordinate transverse to $\ga$. Stretch the neck of $U$ and deform
the equation $\bar{\partial}_{J} u = 0$ to some well-chosen non-homogeneous
equation over the neck in order to force the solutions to approach closed
leaves of $W$. More precisely, for arbitrarily large $\kappa  \ge \de$,
let us denote by $J_{\Si,\kappa}$ any complex structure in the image of
$\si_{norm}$ for the surface $\Si$  which, on $U=\ga \times [-\de, \de]$
is conformally equivalent to the standard complex structure on $\ga \times
[-\kappa, \kappa]$ in such a way that the resulting structure be symmetric
with respect to the $S^1$--action and anti-symmetric with respect to the
$\Z_2$--action that changes the sign of the coordinate in $[-\de, \de]$.
For an arbitrary positive real number $T \ge 0$, let  $h_{T}\co [-\de,
\de] \to \R$ be a non-decreasing function which is constant and equal
to $-T/2$ near $-\de$ and is constant and equal to $T/2$ near the other
end $+\de$; thus its total variation is $T \in \R^+$. Denote by $F_{T}$ its
pull-back to $\pi^{-1}(U)$ by the projections $\pi^{-1}(U) \to \ga \times
[-\de, \de] \to [-\de, \de]$. The generic family $J_b, b \in \Si$,
being fixed, for each  $J_{\Si,\kappa}$, denote by $J_{P,\kappa}$ an
almost-complex structures on $P$, compatible with $\Om$, which restricts
to the family $J_b, b \in \Si$, on the fibers and which is such that
the projection is pseudo-holomorphic with respect with the structures
$J_{\Si,\kappa}$ on the base. For each $\kappa$ and $T$, consider the
sections $u=u_{\kappa,T}\co\Si \to P$ that realise
the class $\tau$, which over
$\ga \times  [-\de, \de]$ are solutions of the equation:
$$\bar{\p}_{J_{\Si,\kappa},J_{\kappa}} u =  - \nabla F_{T}$$
and which over the rest of $\Si$ are solutions of the same equation
and boundary conditions as for the original $u_0=u_{\kappa = \de,
T=0}$. Thus the moduli space of solutions is the set of all quadruples $(
\kappa, J_{\Si, \kappa}, T, u_{\kappa,T})$ where $\kappa \ge \de,
T \ge 0$, $J_{\Si, \kappa}$ belongs to the image of $\si_{norm}$
and $u_{\kappa,T}$ is a solution of the above system. Denote it by
$\Mm_{\kappa}$. Because the solution $u_0$ is both generic and isolated
(amongst {\em all} normalised conformal structures on $\Si$), the
real dimension of $\Mm_{\kappa}$ is $1$. This is because, basically,
the moduli space $\Mm_{\kappa}$ is obtained from the original one by
deforming the right hand side along a one-parameter family given by
$T$. The computation of the
$s$--energy over $\ga \times  [-\kappa, \kappa]$  leads to
$$0 \le \int_{\ga \times  [-\kappa, \kappa]} \left\| \frac{\partial
  u_{\kappa,T}}{\partial s} \right\|^2 \le \left( \int_{\ga \times [-\kappa,
  \kappa]} (u_{\kappa,T})^{*}(\Om) \right)  - T \le \area (u) - T$$
because, by Lemma \ref{le:energy}, the restriction of $\Om$ to the rest of the image of $u_{\kappa,T}$ is
non-negative. Here $\area (u)$ is the $\Om$--area of the solution
$u_{\kappa,T}$ over all of $\Si$, which depends only on the homology
class of the solution and the boundary conditions --  it is a constant
attached to the moduli space. Thus there is no solution when $T > \area
(u)$. By Gromov--Floer's
compactness theorems, the conformal structure of $J_{\Si, \kappa}$  must
degenerate as $T$ approaches its upper bound.  But any degeneracy that
would occur away from the set consisting of the union of the ends and of
$\ga \times  [-\kappa, \kappa]$, is of real codimension $2$ and can
therefore be avoided in $\Mm_{\kappa}$.
Because by our hypothesis on cycles, there is no broken Floer flow line
near the ends (ingoing or outgoing), then the degeneracy must occur
over $\ga \times  [-\kappa, \kappa]$, ie there is a sequence of triples
$(\kappa_i,T_i,u_i)$ such that
$T_i \to T_{\infty} \le \area (u)$ and $\kappa_i \to \infty$.  But because the
$s$--energies of the $u_i$s are bounded, there is a sequence $s_i$
of levels for which the integrals
$$\int_{s=s_i} \left\| \frac{\partial u_i}{\partial s} \right\|^2  dt$$
converge to zero, which means that $J_{\kappa_i} \frac{\partial
u_i}{\partial t} + \nabla F_i$ is arbitrarily $L^2$--small: thus the loops
$u_i |_{s_i}$ must converge to a periodic orbit of the characteristic
foliation of the function $F_{T_{\infty}}\co P \to \R$. But on $U$, this
foliation corresponds to the characteristic flow of a copy of $W$ (lying
over $s_i$). This shows that the space of solutions in the moduli space
$\Mm^{\Ga}_{g,k,\ell,\si_{norm},\{J_b\},F} $ can be deformed to a space
of solutions that has the same counting  and for which each solution
decomposes near $\ga_j$.

This proves the claim if the surgery takes place on a cylindrical end. But
if the surgery takes place on a marked point, one can use exactly
the same ``stretch the neck'' argument with the difference that one
applies it  instead to the symplectically  trivial cylinder $W \times I$
where $I$ is an interval and $W$ is the  trivial product $S^1 \times (M,
\om)$. This produces a {\em flat} holomorphic cylinder over the infinite
cylindrical end, because the closed orbits of a monotone function on
$S^1 \times I$ that depends only on the the variable $t \in I$,  are
the circles $S^1 \times \{t\}$ over which the monodromy in $S^1 \times M
\times I$ is trivial.   By choosing the family $J_b$ over $S^1 \times I$
independent of $b$, and projecting solutions to the fiber, this yields
a solution of an elliptic Cauchy--Riemann type equation over $S^1 \times
[0, \infty)$ and over $S^1 \times (-\infty,0]$ as well. But each one is
conformally equivalent to the unit disc $D^2$ with the origin
removed. Thus, since the solution converges to a constant in the
fiber $M$, one may extend it over all $D^2$ (in other words, we have
reproved here the ``removal of singularity'' theorem for maps with
finite energy). This yields two solutions, one on each unit disc, that
meet transversally at a common point when the two $M$--fibers  over the
origins of the two discs are identified.

(B)\qua It is now easy to generalise the previous argument to the case
when the chains are not assumed to be cycles. Clearly, one can still
consider the same deformed moduli space $\Mm_{\kappa}$.   It is non-empty
and can only degenerate in three ways:
\begin{enumerate}
\item[(1)] at either a Morse or a Floer broken flow line at one of the
ingoing data,
\item[(2)] over some $\ga \times [-\kappa, \kappa]$, or
\item[(3)] at a Morse or a Floer broken flow line at one of the outgoing
data.
\end{enumerate}
But this, counted with signs, leads immediately to
the formula of a homotopy operator:
$$\p H  = H \p + \Phi''_{\si''} \circ  \Phi \circ \Phi'_{\si'} - \Phi_{\si}$$
where the three terms on the right hand side correspond to the three above possibilities in the same order.
\end{proof}

\begin{lemma}
\label{le:identity}
Let $P$ be a topologically trivial $M$--fibration over
$S^2$ with two marked points $p_1 \in \cal A_+$ and $ p_2   \in \cal A_-$,
equipped with a ruled symplectic form $\Om$ deformation equivalent to
the split form, and let
$\si_0$ be the flat section $pt\times S^2$ of $P$.  Then the map
$\Phi_{\si_0}\co QH_*(M) \to QH_*(M)$ corresponding to these marked points
is the identity map.

The same is true if $P$ is a topologically trivial $M$--fibration over
the cylinder,  equipped with a ruled symplectic form $\Om$ deformation
equivalent to
the split form, and if $\si_0$ be the flat section $pt\times (S^1 \times
I)$ of $P$. Then  the map
$\Phi_{\si_0}\co FH_*(M) \to FH_*(M)$ is the identity map.
\end{lemma}

The proof of the first part was sketched in Lemma 4.A of
Lalonde--McDuff--Polterovich \cite{LMP} and was written down for general
symplectic manifolds by McDuff in \cite{McD-sequel}. The second part is
proved similarly -- it is relatively easy in the weakly monotone case
and left to the reader.

 We say that a pair $(P, \si)$ is {\it deformation equivalent to the split
pair} if $P$ has a smooth trivialisation $\phi\co P \to \Si \times M$ that sends
$\si$ to the flat section and the form
$\Om$ to a form which is deformation equivalent to the split form.

\begin{cor}
\label{co:iso}
Let $(M,\om) \hookrightarrow P \to S^2$ be a Hamiltonian fibration with
two marked points $p_1 \in \cal A_+$ and $ p_2 \in \cal A_-$ and $\si$
a homology class of sections up to equivalence. If there is another
pair $(P', \si')$ with two marked points $p'_1 \in \cal A'_+$ and $
p'_2 \in \cal A'_-$ such that the pair $(P \# P', \si \# \si')$ is
deformation equivalent to the split pair, then $\Phi_{P,\si}\co QH(M_{p_1})
\to QH(M_{p_2})$ is an isomorphism.

Let $(M,\om) \hookrightarrow P \to D^2$ be a Hamiltonian fibration with
one marked point $p \in \cal A_+$ and with $\p D^2 \in \cal A_-$ and
let $\si$ a homology class of sections of $\bar P$ up to equivalence. If
there is another similar pair $(P', \si')$ with one marked point $p'
\in \cal A'_-$ and $\p D^2 \in \cal A_+$, such that the monodromies $\al$
and $ \al'$ of their boundaries coincide,  and if the pair $(P\# P', \si
\# \si')$ is deformation equivalent to the split pair, then $\Phi_{\si}\co
QH(M_{p}) \to FH(\p P)$ is an isomorphism.
\end{cor}

\begin{proof}
The proof is a direct consequence of  the Composition theorem
\ref{thm:composition} and of Lemma \ref{le:identity}.
\end{proof}

{\bf Remark}\qua Note that the first part of the corollary says that the
Seidel map is an isomorphism, while the second part implies that the PSS
map from quantum to Floer homology is an isomorphism. Note that in the
second part of the last corollary, one could permute the sets $\cal A_+$
and $\cal A_-$ so that the gluing of $P,P'$ would give a cylinder --
using the second part of Lemma \ref{le:identity} instead, we  would be
led to the conclusion that $\Phi_{\si}\co FH(M_{p}) \to QH(\p P)$ is an
isomorphism, which is then of course the inverse of $\Phi_{\si}\co
QH(M_{p}) \to FH(\p P)$.

\medskip
{\bf Remark}\qua All results of this section obviously  generalise to
ruled symplectic fibrations, ie fibrations $(M, \om) \to P \to \Si$
having $\Symp(M,\om)$ as structural group instead of $\Diff_{Ham}(M,\om)$,
if the following three conditions are satisfied:

\begin{enumerate}
\item[(1)] there is a closed extension of the fiberwise symplectic forms,
\item[(2)] the monodromy round each end is a Hamiltonian diffeomorphism, and
\item[(3)] there is at least one section of the fibration $\bar P$.
\end{enumerate}

Thus, for instance, the monodromy round some closed loop in the base $\Si$
could be a symplectic diffeomorphism not isotopic to the identity. In
this case, the field theory developed in this section cannot be reduced
to a composition of PSS and Seidel homomorphisms.

\section{How essentiality implies minimality}
\label{se:EimpliesM}

The following definition is due to Polterovich \cite{P} and Schwartz
\cite{Sc}:  a generator $(x,v) \in FC_*(H_t)$ is {\em essential} if
there is a class $a \in FH$ such that any cycle $\al $ representing $a$
must contains the element $(x,v)$ (ie it appears with non-vanishing
coefficient in the cycle). An equivalent way of expressing this is:
the inclusion $L \to FC(H_t)$ of the subcomplex generated by all the
generators of $FC_*(H_t)$ except $(x,v)$ induces a morphism between the
corresponding homologies which is not onto.

Let $\Hh$ be the space $C^{\infty}([0,1] \times M, \R)$. It is well-known
that the Hamiltonian paths $\phi^H_{t \in [0,1]}$ generated by these
functions can as well be generated by functions in $C^{\infty}( S^1
\times M, \R)$.

Here and after, we use the following definition of Hofer's length of
an element $H \in \Hh$:  $\int_0^1 (\max_M H_t - \min_M H_t ) dt$. A
function in $\Hh$ is called {\em quasi-autonomous} if there are two
points $P,Q \in M$ such that $P$ is a global maximum of each $H_t, t \in
[0,1]$, and $Q$ is a global minimum of each $H_t, t \in [0,1]$. We will
refer to these points as {\em fixed maximum, fixed minimum} respectively.

For any $H \in \Hh$, denote by $\bar H$ the {\em opposite} Hamiltonian,
ie the function defined by
$$\bar H (t,x) = - H(t, \phi^H_t(x)).$$
This Hamiltonian generates the path $(\phi^H_t )^{-1}$. Obviously, the
roles played by $P$ and $Q$ are reversed while the Hofer length remains
the same.

We will show:

\begin{theorem}
\label{thm:EimpliesM}
Let $(M, \om)$ be a symplectic manifold, that we assume to be weakly
monotone for simplicity. Let $H=H_{t \in [0,1]}$ be any Floer-generic
quasi-autonomous Hamiltonian for which the class $(x_{max}, v_{cst})$ is
essential in $FC_*(H)$ and $(x_{min}, v_{cst})$ is essential in $FC_*(\bar
H)$ (here $v_{cst}$ denotes the constant disc). Then the Hofer length of
the path $\phi_t$ generated by $H$ is minimal among all paths joining the
identity to $\phi_1$ which are homotopic to $\phi_t$ with fixed endpoints
$\id, \phi_1$. If the manifold is symplectically aspherical, the path
$\phi_t$ is length minimising among all paths joining $\id$ to $\phi_1$.
\end{theorem}

Recall that a {\em symplectically aspherical} means that the integral
of $\om$ over any 2--sphere vanishes.

In all results on the minimality of geodesics in Hofer's geometry,
there are essentially two steps: in the first one, one shows that some
conditions concerning the dynamics of the Hamiltonian path (inexistence
of periodic orbits of period one, etc) implies that a more abstract
property is satisfied, like the essentiality of the two classes $(x_{min},
v_{cst})$ and $(x_{max}, v_{cst})$ or some energy-capacity inequality. The
second step then consists in proving that this latter property implies
length-minimality.

The first step was carried out in this setting in \cite{KL} for aspherical
manifolds. The goal of this section is to prove that the second step is
true for all weakly monotone manifolds (we restrict ourselves to weakly
monotone manifolds only for simplicity). There was an attempt  by Oh to
complete this scheme in full generality in \cite{Oh1,Oh2}, but a recent
erratum \cite{Oh-erratum}  restricts the scope of validity of his results
(however the second step in Oh's papers,  by which he proves that the
energy-capacity inequality implies length-minimality, does not seem
at first sight affected by his erratum). In any case,  our aim is to
show that the homological essentiality implies minimality, using only
simple geometric methods and our field theory. The key ingredient of
this section is Proposition~\ref{prop:key}: it is there where our field
theory intervenes. Inasmuch as one is willing to consider that our field
theory can be generalised to all manifolds using perturbations and
virtual cycles, then Theorem~\ref{thm:EimpliesM} would hold for all
symplectic manifolds and therefore the second step of the above scheme
would be completed.

Recall that a fixed point $x$ of the flow $\phi^H_{t \in [0,1]}$ is
{\em under-twisted} if, given any value $T \in [0,1]$, the linearised
flow $D\phi^H_{t \in [0,T]}(x)\co T_xM \to T_xM$ has no non-trivial closed
orbit.  It is  {\em generically under-twisted} if, given any value $T \in
[0,1]$, the linearised flow $D\phi^H_{t \in [0,T]}(x)\co T_xM \to T_xM$
has  only $0$ as fixed point (ie we also exclude all trivial closed
orbits, except the orbit at $0$).

\begin{definition}
Let $H_0$ is a Morse function on a symplectic manifold and let $H_{t \in
[0,1]}$ be a quasi-autonomous Hamiltonian starting with $H_0$. We say
that $H_{t \in [0,1]}$ {\em has no index jump} if
\begin{enumerate}
\item[(1)] each closed orbit at time $1$ is the endpoint of a continuous
family, as $T$ goes from $0$ to $1$, of closed orbits of $\phi^H_{t \in
[0,T]}$ for which the capping discs can be chosen continuously so that,
at time $T=0$, it is the constant one, and
\item[(2)] with respect to these capping discs,  as $T$ goes from $0$ to
$1$ the Conley--Zehnder indices remain the same as in the $t=0$ Morse case.
\end{enumerate}
We say that $H_{t \in [0,1]}$  is {\em restrained} if the condition (1)
holds but the second one is replaced by
\begin{enumerate}
\item[(2$'$)] with respect to these capping discs, the action functional
remains in the range $[\Aa^H_T(Q), \Aa^H_T(P)]$ for all $T \in [0,1]$.
\end{enumerate}
\end{definition}

Note that the hypothesis of  ``under-twisted'' simply means that the
condition (2) holds for the fixed minimum and maximum. Note also that
these conditions do not mean that new closed orbits or index jumps
cannot appear in the time interval $T \in [0,1]$; what they say is that,
if they appear, they must disappear before time $1$.

It is an interesting question to know whether or not each of these two
sets of conditions ((1) and (2), or (1) and (2'))  implies Hofer's
minimality of the path $\phi^H_{t \in [0,1]}$. We will not examine in
this paper the ``restrained'' set of conditions. But we will show that
the first set of conditions (ie (1) and (2)) is sufficient if the
minimal Chern number is either $0$ or large enough:

\begin{cor}
\label{cor:Hofer}
Let $(M^{2n}, \om)$ be a symplectic manifold with minimal Chern number
not in the range $[1, n]$ (a projective space for instance). Suppose that
$H_0$ is a Morse function and let $H_{t \in [0,1]}$ be an undertwisted
quasi-autonomous Hamiltonian  without index jump. Then $H_{t \in [0,1]}$
induces a Hamiltonian path which is minimal in Hofer's length in its
homotopy class with fixed endpoints.
\end{cor}

\begin{proof}
First, it is an easy exercise to check that the proof in Section~8 of
Lalonde--Kerman \cite{KL} applies here as well, with minor changes,
so that we may assume that the Hamiltonian $H_{t \in [0,1]}$ in the
statement of the corollary has the following additional generic
properties:
\begin{enumerate}
\item[(a)] for all $t \in [0,1]$, $P$ is the unique global maximum
of $H_t$ and is non-degenerate,
\item[(b)] for all $t \in [0,1]$, $Q$ is the unique global minimum of
$H_t$ and is non-degenerate, and
\item[(c)] $H$ is Floer-generic.
\end{enumerate}
Now, the generators of the Floer complex $FC_*(H)$
are pairs consisting of a periodic orbit $x(t)$ and a capping disc
$v$. By definition of  ``no index jump'', each closed orbit $x$ is the
endpoint of a continuous family starting at some critical point of $H_0$
with continuously evolving capping discs with constant Conley--Zehnder
index. Thus  a generator $(x,v)$ has index equal to the Morse index
of a critical point of $H_0$ (in $[0,2n]$) plus or minus some multiple
of twice the Chern class of some 2--sphere in $M$. By hypothesis, this
Chern number cannot be in $[1,n]$, so the index of a generator cannot
be equal to $2n+1$. This shows that two different cycles of index $2n$
cannot be homologous. Hence, to establish the essentiality of $(P,
v_{cst})$ in $FC_*(H)$ (here $v_{cst}$ denotes the constant disc at
$P$), there only remains to show that there is a cycle in $FC_{2n}(H)$
that contains the element $(P,v_{cst})$.  This is done exactly as in
the proof of Proposition~5.2 of \cite{KL}. (It is for that proof that one
needs the above additional generic properties (a) and (b).) Actually,
the proof in \cite{KL} is written down for an aspherical manifold, but
it works as well for a general manifold; one simply needs to replace $P$
by $(P,v_{cst})$. This establishes the essentiality of $(P, v_{cst})$
in $FC_*(H)$; the essentiality of $(Q,v_{cst})$ in $FC_*(\bar H)$ is
proved similarly.
By Theorem~\ref{thm:EimpliesM} above, we get minimality.
\end{proof}

Here is a simple example of application of that Corollary.  Given a
symplectic manifold with minimal Chern number not in the range $[1, n]$,
$\CP^n$ for instance, let $f\co M \to \R$ be a Morse function, for instance
a function on $\CP^n$ which once pulled-back on $\C^{n+1}-\{0\}$, is of
the form
$$\frac{\sum_{i,j} (a_{i,j} x_ix_j + b_{i,j} x_i y_j + c_{i,j}
  y_iy_j)}{\|z\|^2}.$$
Let $P$ be a global maximum of $f$ and $Q$ a global
minimum of $f$ such that the linearised flows have no nontrivial closed
orbit in time $\le 1$.  Note that they are then automatically generically
undertwisted by the Morse condition. Let $\eps > 0$ be sufficiently small
so that the flow $\phi^f_t$ induced by $f$ be transversal to the diagonal
in $M \times M$ for all $t \in (0,\eps]$. Consider any Hamiltonian $H_{t
\in [0,1]}$ such that:
\begin{enumerate}
\item[(1)] on some open neighbourhoods of  $P$ and $Q$, $H_t = f$ for
  all $t \in [0,1]$,
\item[(2)] $H_t = f$ for $t \le \eps$, 
\item[(3)] the graph of $\phi^H_t$ remains transversal to the diagonal
  in $M \times M$ for all $t \in [\eps,1]$, and
\item[(4)] $H_t(x) \in [f(Q), f(P)]$ for all $x \in M$ and $t \in [0,1]$.
\end{enumerate}

Then, clearly, each fixed point of $\phi^H_1$  is the endpoint of a
continuous family, as $T$ goes from $0$ to $1$, of closed orbits of
$\phi^H_{t \in [0,T]}$. Thus the capping discs can be chosen continuously
so that, at time $T=0$, they are the constant ones. By the preceding
corollary, one concludes that this path is length minimising rel endpoints
in its homotopy class.

Let us go back to Theorem~\ref{thm:EimpliesM}.
Its proof takes the next two paragraphs: in the first one,  we introduce
the geometric tools that will be necessary:  gluing of monodromies,
Hamiltonian cylinders and fibrations, and area estimates. In the
second paragraph, we  conclude the argument in two cases:
\begin{enumerate}
\item[(1)] $M$ is monotone and the paths are in the same homotopy class
rel endpoints in $\Ham(M,\om)$;
\item[(2)] $M$ is aspherical
\end{enumerate}
and we prove the stronger statement
of Theorem~\ref{thm:EimpliesM}, ie minimality with respect to {\em all}
paths with fixed endpoints.

\subsection{Gluing Hamiltonian fibrations along monodromies}
\label{ss:monodromies} \label{sse:hofer}

Let us first recall from \cite{LM-Invent, LM-GAFA} that if $H_{t \in
[0,1]}, K_{t \in [0,1]}$ generate paths in the group of Hamiltonian
diffeomorphisms of $(M, \om)$ joining the identity to the same Hamiltonian
diffeomorphism $\phi = \phi^H_{t=1} = \phi^K_{t=1}$, one can construct a
symplectic manifold  $R_{H,K}$  by gluing the region under the graph of
$H$ with the region above the graph of $K$. We will introduce a variant
of this construction.

First, after reparametrisation of the Hamiltonian paths $\phi_{t \in
[0,1}$, we may assume that all the generating  Hamiltonians functions
$G_{t \in [0,1]}$ are normalised so that
\begin{enumerate}
\item[(1)] they are the restrictions to the time interval $[0,1]$ of a
  time-periodic Hamiltonian of period $1$ and
\item[(2)] the minimum over $M$ of $G_t$ equals $0$ for every  $t$
  and the maximum over $M$ of $G_t$ is independent of $t$ (it is therefore
  equal to $\Ll(G)$).
\end{enumerate}
Then $G$ may be considered as a map
$G\co M \times S^1 \to [0,\infty)$.  The restriction of the form
$\Om = \om \oplus ds \wedge dt \in \Om^2(M \times S^1 \times \R)$ to the graph
of a normalised Hamiltonian $G$ gives rise to a characteristic foliation. It is
well known that the monodromy of that foliation is equal to
$\phi^G_1$. The ``region over the graph'' of
$G$ is the subset of $M \times S^1 \times \R$ defined by
$$R^+(G) = \left\{(x,t,s) : G_t(x) \le s \le \max_M G_t = \Ll(G) \right\}$$
and the ``region under the graph'' of $G$ is the subset
of $M \times S^1 \times \R$
$$R^-(G) = \left\{(x,t,s) : \min_M G_t = 0 \le s \le G_t(x) \right\}.$$
In order that all constructions in the sequel be smooth, we will
have to consider spaces $R^-_{\eps}(G)$ for
arbitrarily small $\eps > 0$:
$$R^-_{\eps}(G) = \{(x,t,s) : 0 \le s \le G_t(x) + \eps \}.$$
so that the lower and upper components of the boundary of $R^-_{\eps}(G)$ are disjoint.
Similarly:
$$R^+_{\eps}(G) = \{(x,t,s) : G_t(x) \le s \le \Ll(G) + \eps\}.$$
Note that $R^-_{\eps}(G)$ is symplectomorphic to $\{(x,t,s) : -\eps \le s
\le G_t(x) \}$, so both spaces are obtained by gluing to $R^-(G),R^+(G)$ the
product of $(M,\om)$ with the annulus
$S^1 \times I$, where $I$ is an interval of length $\eps$.

Let $q_-\co M \times S^1 \times [0,\infty) \to M \times \R^2$
be the map that is the identity on $M$, sends $t$ to the angle coordinate of
$\R^2$ that we will still denote by $t$, and maps
$s \in [0,\infty)$ to the action coordinate $c = \pi r^2
\in [0,\infty)$ of $\R^2$. Thus $q_-$ maps the form $\om + ds \wedge dt$ to the
form $\om + dc \wedge dt$  and sends $R^-_{\eps}(G)$
to
$$P^-_{\eps}(G) = \{(x,t,c) : 0 \le c \le G_t(x) + \eps \} \subset M
\times \R^2.$$
Similarly, let $q_+\co M \times S^1 \times (-\infty, \Ll(G) + \eps]
\to M \times \R^2$ be the map that is the identity on $M$, sends $t$
to the angle coordinate of $\R^2$ and maps
$s$ to $c(s)=\Ll(G) + \eps -s$. Thus $q_+$ maps the form $\om + ds \wedge dt$ to
$\om - dc \wedge dt$  and sends $R^+_{\eps}(G)$ to
$$P^+_{\eps}(G) = \{(x,t,c) : 0 \le c \le \Ll(G) + \eps - G_t(x) \}
  \subset M \times \R^2.$$
This means that, with respect to the orientation ${\cal O^+}$ induced by $- dc
\wedge dt$ on the boundary of the base of $P^+_{\eps}(G)$, the monodromy of the
characteristic flow of $\p P^+_{\eps}(G)$, followed in the direction of ${\cal
O^+}$, is the inverse of the monodromy of $\p P^-_{\eps}(G)$, followed in the
direction  ${\cal O^-}$ induced by $dc \wedge dt$ on the boundary of the base of
$P^-_{\eps}(G)$.

Hence, up to diffeomorphism, $P^{\pm}_{\eps}(G)$ may be considered as
a topologically trivial fibration $M \stackrel{\pi}{\hookrightarrow}
P^{\pm}_{\eps}(G) \to D^2$ equipped with a ruled symplectic form $\Om$,
ie a symplectic form that restricts to a non-degenerate form on each
fiber, such that the monodromy round the boundary of the base covered
in the positive direction is $\phi^G_1$ for $P^-_{\eps}(G)$ and
$(\phi^G_1)^{-1}$ for $P^+_{\eps}(G)$.
Here the positive direction is intrinsically defined by following the symplectic
gradient of a function whose level set is $\p P^{\pm}_{\eps}(G)$ and whose
ordinary gradient points outward. The diffeomorphism that does this will
be denoted by $\al^-_{\eps}(G)$: it preserves $t$
and maps the sets $(x,t, G_t(x) + \eps)$, for each given $t$, to the $M$
fiber over the point of $\partial D^2$ with angle $t$ (and similarly
for $\al^+_{\eps}(G)$).

Note finally that any ruled symplectic form $\Om$ on a fibration $M
\hookrightarrow P \to D^2$ has a fiberwise symplectic trivialisation, ie a map
$\phi\co P \to D^2 \times (M, \om)$ that commutes with the projection
to $D^2$ and maps $\Om$ to a form that restricts  precisely to $\om$
on each fiber (this trivialisation can be constructed by using the
connexion induced by $\Om$ over the rays in $D^2$ emanating from a base
point). However, if $\pi_1(\Ham(M,\om))$ is not $\{0\}$, this
trivialisation is not unique up to homotopy. Our construction associates
to each normalised Hamiltonian $G$ a ruled symplectic fibration
$P^{\pm}_{\eps}(G)$ equipped with such a trivialisation.

Let $P_1,P_2$ be two $M$--fibrations over $D^2$ equipped with ruled
symplectic structures $\Om_i$ and with the fibers over the base point $1 \in
\p D^2$  symplectically identified. Suppose that the monodromies, via this
identification, are inverses of each other. This yields a ruled
symplectic structure on a $M$--fibration $P_1 \# P_2$ over
$S^2$ by gluing the monodromies in the obvious way. This applies in particular
to the pair $P^-_{\eps}(H), P^+_{\eps}(K)$ if $H_t, K_t$ generate Hamiltonian
flows with the same time-one maps $\phi^H_{t=1} = \phi^K_{t=1}$. Denote by
$P_{\eps}(H,K)$ the resulting fibration over $S^2$. Similarly, denote by
$P_{\eps}(K,H)$ the fibration over $S^2$ obtained by gluing $P^-_{\eps}(K)$ to
$P^+_{\eps}(H)$. We will call them the {\em mixed fibrations} associated to $H$
and $K$. Note that $P_{\eps}(G,G)$ is the symplectically trivial fibration over
a sphere of area $\Ll(G) + 2 \eps$.

\begin{definition}\qua
\label{def:P}

(i)\qua Let $(M,\om) \hookrightarrow  (P, \Om) \to B$ be a
ruled symplectic form over a compact surface $\Si$. The {\em area} of $P$
is by definition the quotient of the $\Om$--volume of $P$ by the $\om$--volume
of the fiber.

(ii)\qua When $M \to P \to D^2$ is equipped with a fiberwise symplectic
trivialisation,   the characteristic flow round $\p P$ gives rise to a
Hamiltonian flow $\phi_t$ on $M$. Then each pair $(x(t), v)$ of a periodic
orbit of $\phi_t$ and of a disc $v$ of $M$ (defined up to homotopy)
with boundary equal to $x(t)$ corresponds via the trivialisation to
a section of $P$, that we will denote $\si(x(t),v)$, with boundary
lying on the characteristic leaf of $\p P$ given by $x(t)$.   If $P_1,P_2$
are two fibrations with fiberwise symplectic trivialisations and
inverse monodromies, a periodic orbit $x(t)$ of the flow $\phi^{P_1}_t$,
a bounding disc $v_1 \subset M$, and a bounding disc $v_2 \subset M$ of
the flow $\phi^{P_2}_t (x(0))$ gives rise to a section $\si(x(t),v_1,v_2)$
of $P_1 \# P_2$.
\end{definition}

{\bf Remark}\qua Note that if
$H_{t \in [0,1]}$ has a fixed maximum $p_{max}$ and if $v_{p_{max}}$
denotes the constant disc, the Hofer length of $H_{t \in [0,1]}$
is of course equal, up to $\eps$, to the $\Om$--area of the  section
$\si(p_{max}, v_{p_{max}})$ in $P^-_{\eps}(H)$. The same applies for
$p_{min}$ and $P^+_{\eps}(H)$.

\begin{prop}[Basic inequalities]
\label{prop:basic}
Let $H_t, K_t$ generate Hamiltonian flows with the same time-one maps
$\phi^H_{t=1} = \phi^K_{t=1}$.  Then $\Ll(H_t) \le \Ll(K_t)$ if and only if
$$2 \Ll(H_t)  \le \area (P_{\eps}(H,K)) + \area (P_{\eps}(K,H))$$
for all $\eps > 0$. Therefore $\Ll(H_t) \le \Ll(K_t)$ if each
of the following inequalities holds for all $\eps > 0$:
\begin{eqnarray}
\label{eqn:1}
\Ll(H_t) \le \area (P_{\eps}(H,K)) \\
\label{eqn:2}
\Ll(H_t) \le \area (P_{\eps}(K,H))
\end{eqnarray}
Moreover, if the space $P_{\eps}(H,K)$ has a homology class
of sections
$\si_{flat}^{H,K}$ whose
$\Om_{H,K}$--area is equal to $\area(P_{\eps}(H,K))$ and which decomposes as
$$\si_{flat}^{H,K} = \si(p_{max},v_{p_{max}}) + \si'$$
over $P^-_{\eps}(H)$ and $P^+_{\eps}(K)$ respectively, then the inequality
\eqref{eqn:1} amounts to
$$\mbox{$\Om_{P^+(K)}$--area of $(\si') \;  \ge \; 0$}.$$
Similarly, if the space $P_{\eps}(K,H)$ has a homology class of sections
$\si_{flat}^{K,H}$ whose
$\Om_{K,H}$--area is equal to $\area(P_{\eps}(K,H))$ and which decomposes as
$$\si_{flat}^{K,H} =   \si' + \si(p_{min},v_{p_{min}})$$
over $P^-_{\eps}(K)$ and $P^+_{\eps}(H)$ respectively, then the inequality
\eqref{eqn:2} amounts to
$$\mbox{$\Om_{P^-(K)}$--area of $(\si') \; \ge \;  0$}.$$
\end{prop}

In that statement, of course, the various indices affecting $\Om$
denote the spaces in which the form lives (we deleted the $\eps$ in
the indices).

\begin{proof}
We have:
$\Ll(H_t) \le \Ll(K_t)$ if and only if
$$2 \vol (P_{\eps}(H,H)) \le \vol (P_{\eps}(K,K)) + \vol (P_{\eps}(H,H)).$$
But the latter is equal to $\vol(P_{\eps}(H,K)) + \vol
(P_{\eps}(K,H))$, and after dividing out by the volume of $M$, we find $2 \Ll
(H_t) + 4 \eps
\le \area (P_{\eps}(H,K)) + \area (P_{\eps}(K,H))$ for all $\eps > 0$, which
means that
$$2 \Ll(H_t) \le \area (P_{\eps}(H,K)) + \area (P_{\eps}(K,H))$$
for all $\eps > 0$.

The rest of the proposition is a direct
consequence of the remark preceding the proposition.
\end{proof}

This establishes the strategy for the proof of minimality: we will find
sections of the ruled symplectic fibrations $P_{\eps}(H,K), P_{\eps}(K,H)$
which decompose as above and for which the $\si'$--part is positive.

We need a last result on the relation of the areas of the fibrations
$P_{\eps}(H,K))$, $P_{\eps}(K,H)$ and the areas of {\em flat sections} in
these spaces. We will consider two cases:
\begin{enumerate}
\item[(1)] the path $\phi^K_t$ is homotopic to $\phi^H_t$, and
\item[(2)] the condition on the homotopy of the paths is released but $M$
  is aspherical.
\end{enumerate}

Let $(P, \Om)$ be a $(M, \om)$--ruled symplectic manifold over the
$2$--sphere.  By the characterisation of Hamiltonian bundles in \cite{MS},
Theorem~6.36 (or \cite{LM-Topology}, Theorem~1.1), any ruled symplectic
manifold over a simply connected base is a Hamiltonian fibration. Now
suppose that there is a fiberwise symplectic trivialisation $\psi\co
P \to S^2 \times M$.  Then any other such trivialisation $\psi'$ will
differ by a map $\Psi\co S^2 \to \Symp_0(M)$. Because $\pi_2(\Symp_0(M))
= \pi_2(\Ham(M))$, $\Psi$ is homotopic to a map $S^2 \to \Ham(M)$. But
these act trivially on the homology of $S^2 \times M$ by Proposition~5.4
in Lalonde--McDuff \cite{LM-Topology}. This means that the homology class
of the section $\psi^{-1} (S^2 \times \{pt\})$ is independent of the
chosen trivialisation of $P$. Such a section will be called {\em flat}.

Now, if $\phi^K$ is homotopic to $\phi^H$ with fixed endpoints, the
fibrations $P_{\eps}(H,K)$ and $P_{\eps}(K,H)$ are symplectically
trivial. The fiberwise symplectic trivialisations are  induced by
a choice of a homotopy $G_{s,t}$ between the Hamiltonian paths $H_t$
and $K_t$. Thus there is a well-defined class of flat sections. Here is
a description of that class: if $G_{s,t}$ is a homotopy of paths with
fixed endpoints between $H$ and $K$, with $G_{0,t} = H, G_{1,t} = K$,
and if $(x(t),v_x)$ is any pair of a closed orbit of $H$ and a bounding
disc, then the homotopy $G_{s,t}$ gives rise to a homotopy between
$(x,v_x)$ and $(\phi^K_t(x), v^K_x)$, which defines a bounding disc
$v^K_x$ up to homotopy (ie as an element in $\pi_2(M,\phi^K_t(x))$).
Both $(x,v_x)$ and $(\phi^K_t(x), v^K_x)$ may be interpreted as sections
of the corresponding cylinders $P^-_{\eps}(H), P^+_{\eps}(K)$, which
after gluing, give the flat section $\si_{flat}$ of $P$.

\begin{lemma}
Let $(M,\om) \hookrightarrow (P,\Om) \to S^2$ be a  ruled symplectic
manifold with compact fiber, with $P = P_{\eps}(H,K)$ or $P_{\eps}(K,H)$
where $K$ is a path homotopic to $H$ with fixed endpoints, and $H$ has
a fixed minimum and fixed maximum.  Then $P$ is a Hamiltonian fibration
and there is a unique homology class $\si_{flat} \in H^2(P)$ which
represents the flat section.  In the case of $P_{\eps}(H,K)$, that class
can be decomposed as the union (along the common boundary loop) of two
homologically well-defined sections of the spaces $R^-(H)$ and $R^+(K)$,
the first one given by $(p_{max},v_{p_{max}})$, with  $v_{p_{max}}$
the constant disc, and the second one given by the corresponding pair
$(\phi^K_t(p_{max}), v^K_{p_{max}})$.  Moreover the area of $P$ is equal
to the $\Om$--area of $\si_{flat}$:
$$\vol (P, \Om) = (\vol (M,\om)) \int_{\si_{flat}} \Om.$$
A similar statement applies to $P_{\eps}(K,H)$ with $p_{max}$ replaced by
$p_{min}$.
\end{lemma}

This lemma produces the $\si'$ required in the last proposition, and
there will remain to show later (see Section~\ref{ss:end}) that its area
is non-negative.

\begin{proof}
There remains only to establish the last equality. Here is a proof using
Moser's argument. Take any homotopy $G_{s,t}$ and consider the
corresponding spaces $P_s, \Om_s$ obtained by gluing the cylinders
$P^{\pm}_{\eps}(G_0) = P^{\pm}_{\eps}(H)$ and $P^{\mp}_{\eps}(G_s)$. This
gives a one-parameter family of forms $\Om_s$ on $S^2 \times M$ for
which some fixed fiber, say $M_0$, is symplectic for all $0 \le s \le
1$. Since $ \vol (P, \Om_s) - (\vol (M,\om)) \int_{\si_{flat}} \Om_s$ is a
continuous function of $s$, the inverse image of $0$ is closed and
therefore, if this equality fails to be true for all $s$, there is smallest
$s_0$ for which it does not hold. But, at this value $s_0$, one may
transform the deformation $\Om_s$ into a genuine isotopy on some small
interval $[s_0, s_0 + \eps)$ by either adding or subtracting a small
multiple of the Thom class of the normal bundle of the fiber $M_0$. This
boils down to add or remove a small symplectically split tube $D^2
\times M_0$. These new forms are isotopic to $\Om_{s_0}$ and therefore
the equality holds. But, it holds for the small (added or removed)
split tube too, so it must hold for the undeformed form as well.
\end{proof}

Let us consider now the second case, when $M$ is aspherical but the path
$\phi^K_t$ need not be homotopic to $\phi^H_t$. The mixed fibration $P =
P_{\eps}(H,K)$ or $P_{\eps}(K,H)$ has therefore a single homology class
of sections. The unicity is clear from the aspherical condition. The
existence is a consequence of the fact that in this case $P$ is a
Hamiltonian fibration corresponding to the loop $(\phi_t^K)^{-1} \star
\phi_t^H  \in \pi_1(\Ham(M))$ or its inverse, and it is a consequence
of the Arnold conjecture that the orbit of such a loop on a point of $M$
is contractible (see \cite{LMP} for instance).

\begin{prop}
Let $(M,\om)$ be an aspherical manifold and $H_t, K_t$ two
Hamiltonians generating paths with the same endpoints, with $H$
quasi-autonomous.  Then each of the spaces
$P = P_{\eps}(H,K)$ or $P_{\eps}(K,H)$ has a single homology class of sections.
Moreover, denoting by $\si_{H,K}, \si_{K,H}$ these sections, we have:
$$\area (P_{\eps}(H,K)) + \area (P_{\eps}(K,H)) = \int_{\si_{H,K}} \Om_{H,K} +
  \int_{\si_{K,H}} \Om_{K,H}.$$
\end{prop}

\begin{proof}
Consider an arc $\ga \subset M$ joining a fixed minimum $p_{min}$ of
$H_t$ to a fixed maximum $p_{max}$ of $H_t$. Its trace under the flow
of $\phi^H$ defines a subset of the graph of $H$ (as a subset of $M
\times \R^2$, with action-angle $c,t$--coordinates):
$$\begin{array}{ccc}
f\co [0,1] \times [0,1] & \longrightarrow & M \times \R^2  \\
 (z, t)              & \longmapsto &
  \left(\phi^H_t\left(\ga\left(z\right)\right), t,
    H_t\left(\phi^H_t\left(\ga\left(z\right)\right)\right)\right)
\end{array}$$
Because $M$ is aspherical, there is up to homotopy a unique way to
extend this map to a cylinder
$$C_H\co [0,1] \times S^1 \to M \times \R^2$$
which is obtained by gluing $f$ with a map
$$f'\co [0,1] \times [0,1]  \to  M \times \{(c,t) : t=0\}$$
in such a way that the projection of $C_H$ on the first factor is a
homotopically trivial $2$--sphere. Note that, when considered as a subset
of $P_{\eps}(H,H) = M \times S^2$, this cylinder can be extended in the
obvious way to a 2--sphere in the flat section class by capping small
discs of area $\eps$.The $\Om_{H,H}$--area of $C_H$ is the length of
$H_{t \in [0,1]}$. Because the graph of $H$ is mapped symplectically
to the graph of $K$ along the characteristic foliations via
the diffeomorphisms $\al^{\pm}_{\eps}(H), \al^{\pm}_{\eps}(K)$
(see the construction described in the three paragraphs before
Definition~\ref{def:P}), the cylinder $C_H$ is mapped to a cylinder $C$
with same area in the space $\graph (K) \subset P_{\eps}(K,K)$ . Now $C$
can be extended in a unique way, up to homotopy, to a 2--sphere $C_K$
of $P_{\eps}(K,K) = (M \times S^2, \Om=\om \oplus \si)$ where $\si$
is the standard area form on $S^2$ whose area is the  $\Om$--area of the
fibration ie is $\Ll(K) + 2\eps$. It is obtained by gluing two discs to $C$:
$$g^+, g^-\co D^2 \to P^+_{\eps}(K), P^-_{\eps}(K).$$
Here $g^+$ is uniquely determined up to homotopy by requiring that its
boundary  is mapped to $(\phi^K_t(p_{max}), t, K_t(\phi^K_t(p_{max})))
\subset \graph (K)$.
Similarly, $g^-$ has its image in $P^-_{\eps}(K)$ and its boundary is sent
to $(\phi^K_t(p_{min}), t, K_t(\phi^K_t(p_{min})) ) \subset \graph (K)$.

This construction shows that (we omit the epsilons for simplicity)
\begin{multline*}
\Ll(K) = \area (C_K) = \area (C) + \area (g^+) + \area (g^-) \\
  = \Ll(H) + \area (g^+) + \area (g^-).
\end{multline*}
   On the other hand, the decomposition of the flat section of $P_{\eps}(H,K)$ described in the last lemma leads to:
$$\int_{\si_{H,K}} \Om_{H,K} =  \Ll(H) + \int_{g^+} \Om_{H,K}$$
and the decomposition of the flat section of $P_{\eps}(K,H)$ to:
$$\int_{\si_{K,H}} \Om_{K,H} =  \int_{g^-} \Om_{H,K} + \Ll(H).$$

Putting this together, we have:
\begin{multline*}
\int_{\si_{H,K}} \Om_{H,K} + \int_{\si_{K,H}} \Om_{K,H} =
  \Ll(H) + \left(\Ll(H) + \int_{g^+} \Om_{H,K} + \int_{g^-} \Om_{H,K}\right) \\
= \Ll(H) + \Ll(K) = \frac{\vol (P_{H,H}) + \vol (P_{K,K})}{\vol (M)} \\
= \frac{\vol (P_{H,K}) + \vol (P_{K,H})}{\vol (M)}
  = \area (P_{H,K}) + \area (P_{K,H}).
\end{multline*}
\end{proof}
As a consequence of that proposition and Proposition~\ref{prop:basic}, we then
have:
\begin{cor}
\label{cor:aspherical}
In the aspherical case, the path
$\phi^H_{t \in [0,1]}$ is minimal amongst all paths with same endpoints if
$$\int_{\si_{H,K}} \Om_{H,K} \ge \Ll(H) \quad {\mbox and} \quad
\int_{\si_{K,H}} \Om_{K,H} \ge \Ll(H)$$
for all $\eps > 0$.
\end{cor}

\subsection{End of the proof} \label{ss:end}

\begin{prop}
If $H_{t \in [0,1]}, K_{t \in [0,1]}$ induce paths with same
endpoints, $\si$ is any section class of $P= P_{H,K}$, and $p_1 \in \cal
A_+$ and $p_2 \in \cal A_-$ are two marked points in $P$, then
$$\Phi_{P,\si}\co QH(M_{p_1}) \longrightarrow QH(M_{p_2})$$
is an isomorphism.
The same statement holds for $P = P_{K,H}$.
\end{prop}

\begin{proof}
This is a consequence of Corollary \ref{co:iso}. Indeed, compose $P_{H,K}$
with $P_{K,H}$. We then get, inside $M \times S^2$, the graphs of two
Hamiltonian paths that are the inverse of each other. Isotope the union
of the two graphs to the graph of the trivial loop in $\Ham(M,\om)$.
Finally choose the section of $P_{H,K}$ so that its connected union with
$\si$ is the flat section.
\end{proof}

The next proposition is the key result of this section:

\begin{prop}
\label{prop:key}
Let $H_t,K_t$ be two Hamiltonian paths with same endpoints on a symplectic
manifold. Assume either that they are homotopic rel endpoints or that
$M$ is aspherical. If a pair $(x(t), v_x)$ is essential in the Floer
homology of $H_t$, then there is a solution $u$ in $P^{\pm}_{\eps}(K)$
to the equation $(**)$, which converges to $\phi^K_t(x(0))$ in the obvious
class $v_K$.
\end{prop}

\begin{proof}
To fix notations, consider $P^{+}_{\eps}(K)$. If the statement in the
proposition did not hold, the isomorphism
$$\Phi_{\si}\co QH_*(M_{p_1}) \to QH_*(M_{p_2})$$
in the bundle $P = P_{\eps}(H,K)$ over $S^2$ would factorise by the Composition Theorem \ref{thm:composition} through 
$$QH_*(M_{p_1}) \stackrel{\Phi'_{\si'_{flat}}}{\longrightarrow}
  FH_*(\p P^{-}_{\eps}(H)) \stackrel{\Phi}{\longrightarrow}
  FH_*(\p P^{+}_{\eps}(K)) \stackrel{\Phi''_{\si''_{flat}}}{\longrightarrow}
  QH_*(M_{p_2})$$ 
where it would be enough to consider a subcomplex
$(L, \p)$ of $FC_*(\p P^{-}_{\eps}(H) \subset P^{-}_{\eps}(H))$ that does not
contain $(x(t), v_x)$. Thus $\Phi_{\si}$ would factorise through $H_*(L)$. Since
both the domain and codomain of $\Phi_{\si}$ can be canonically identified with
$HF_*(H_t)$, this would mean that the inclusion of $L$ in $FC_*(H_t)$ would
induce an isomorphism, a contradiction with the hypothesis of essentiality.
\end{proof}

   We proved the last proposition using the Composition Theorem. It is clear that, unwrapping the proof of the Composition Theorem,  the last proposition is actually established by a ``stretch-the-neck'' argument.

\medskip
{\bf End of the proof of Theorem~\ref{thm:EimpliesM}}\qua
We will apply the previous proposition to the pairs $(x_{min},v_{cst})$
and $(x_{max},v_{cst})$ (where $v_{cst}$ is the constant map), which
are essential in the Floer homology of $H_t$. By
Proposition~\ref{prop:basic} and Corollary~\ref{cor:aspherical}, it is
enough to show:
\begin{enumerate}
\item[(1)] the symplectic area of the solution $u_{max}$ of the
preceding proposition to the equation $(**)$ in  $P^{+}_{\eps}(K)$ that
converges to $\phi^K_t(x_{max})$ in class $v_K$ is non-negative, and
\item[(2)] the symplectic area of the solution $u_{min}$ of the preceding
proposition to the equation $(**)$ in  $P^{-}_{\eps}(K)$ that converges to
$\phi^K_t(x_{min})$ in class $v_K$ is non-negative.
\end{enumerate}
  To prove (1), consider the
$s$--energy of $u_{max}$:
\begin{multline*}
0 \le \int_{\Si} \left\| \frac{\p u_{max}}{\p s} \right\|^2
  = \int_{\Si} \left\langle \frac{\p u_{max}}{\p s},
    -J \frac{\partial u}{\partial t}  - \nabla F \right\rangle \\
\le \Om_{P^{+}_{\eps}(K)} \, \area (u_{max}) - \Totvar (f)
\end{multline*}
Here, on the end (identified with the semi-infinite cylinder),
$\frac{\p}{\p s}$ and $\frac{\p}{\p t}$ are the standard vector fields
while, inside the unit disc, they are given as $r dr$ and $d \theta$
respectively. Because ${\rm Totvar} (f)$ is non-negative, we get:
$$\Om_{P^{+}_{\eps}(K)} \, \area (u_{max}) \ge 0.$$
A similar argument applies to $(x_{min},v_{cst})$.
\endproof

\section{Higher-genus norms}

As an addendum to this paper, we briefly discuss a semi-norm of higher
genus on the group of Hamiltonian diffeomorphisms and state a conjecture
on its relation with the Hofer norm.

Let $(M,\om) \hookrightarrow P \to \Si_{g,1}$ be a Hamiltonian fibration
over a compact  oriented surface obtained by removing an open disc from
a closed surface of genus $g$. Let $q$ be a base point on $S = \p \Si$,
and $\phi \in \Diff_{\Ham}(M,\om)$ the monodromy of the characteristic
foliation on $W=\p P$.

\begin{definition}
Let $\phi \in \Diff_{\Ham}(M,\om)$ be given. Its {\em genus $g$ norm}, $\|
\phi \|_g $ is by definition the infimum, over all Hamiltonian fibrations
$(M, \om) \hookrightarrow P \to \Si_{g,1}$ with monodromy equal to $\phi$,
of the area of $P$.
\end{definition}

Let ${\cal R}$ be the semigroup of all non-increasing sequences of
non-negative real numbers, which contain only finitely non-zero terms,
with the operation:
$$(a_0, a_1, a_2, \ldots )  +_T  (b_0, b_1, b_2, \ldots)
  = (c_0, c_1, c_2, \ldots)$$
where the $c_k$s are defined by: 
$$c_k = \min_{0 \le j \le k} (a_j + b_{k-j}).$$
This operation is associative and commutative. Because the sequences
are non-increasing, the zero sequence is a neutral element. Defining
the {\it total norm} of a Hamiltonian diffeomorphism as
$ \| \phi \|_T = (\| \phi \|_0, \| \phi \|_1, \| \phi \|_2, \ldots )$,
one can show easily:

\begin{prop}
The total norm $\| \cdot \|_T \in {\cal R }$ satisfies the triangle
inequality with respect to the total sum  defined above, ie
$$\| \phi \circ \psi \|_T \le \| \phi \|_T +_T  \| \psi \|_T.$$
\end{prop}

\begin{prop}
\label{prop:normg}
Let $\phi \in \Diff_{\Ham}(M,\om)$, different from the identity. Then
$$\| \phi \|_g = 0$$
if and only if $g \le \ell_c$.
\end{prop}

\begin{proof}
To simplify notation, we will first prove this when $\ell_c = 1$. In
this case, $\phi = [f,g]$.
Endow $([0,\eps] \times [0,\eps]) \times M$ with the split symplectic
structure, identify
$[0, \eps] \times \{0\} \times M$ to $[0, \eps] \times \{\eps\} \times M$
via the map $f$ to get a Hamiltonian fibration over a cylinder. Identify
$\{0\} \times [\eps/3, 2\eps/3] \times M$ with
$\{\eps\} \times [\eps/3, 2\eps/3] \times M$ via the map $g$. This
gives a Hamiltonian  fibration over the punctured $2$--torus whose
monodromy is $[f,g]$. Its area is arbitrarily small, thus $\| \phi
\|_1 = 0$. The same argument shows that if $\phi$  can be written as a
product of $k$ commutators, one can construct a Hamiltonian fibration
$P \to \Si_{k,1}$ of arbitrarily small area whose monodromy round the
boundary is $\phi$. Thus $\| \phi \|_k =0$. Now, to prove that all other
semi-norms $\| \phi \|_{k'}$ also vanish for $k' > k$, one can add as
many small handles as one wants over an arbitrarily small disc where
the above fibration $P \to \Si_{k,1}$ is trivialised.

The converse statement saying that $\| \phi \|_k \neq 0$ if $\phi$ has
commutator length larger than $k$ is a consequence of the non-degeneracy
of the usual Hofer norm and will be proved elsewhere (see \cite{LAR}).
\end{proof}

 Because $\Diff_{\Ham}(M,\om)$ is a simple group, the group generated by products of commutators must be the whole of $\Diff_{\Ham}(M,\om)$. Thus any element $\phi \in \Diff_{\Ham}(M,\om)$ can be written in the form:
$$\phi = [\phi_1,\phi'_1] \circ \ldots \circ  [\phi_k,\phi'_k].$$
\begin{definition}
The least such integer $k$ is called the {\em commutator length} of
$\phi$, $\ell_c(\phi)$.
\end{definition}

Let $\Diff_{\Ham}^k(M)$ be the subspace of all diffeomorphisms with
commutator length less or equal to $k$. The following conjecture is
studied in \cite{LAR}.

\begin{conj}
Let $\phi \in \Diff_{\Ham}(M)$. Then the distance, in the usual Hofer
norm, between $\phi$ and $\Diff_{\Ham}^k(M)$ is equal to $\| \phi \|_k$.
\end{conj}

The last Proposition and Conjecture are closely related to Entov's
paper on the relation between $K$--area and commutator length. Indeed,
Proposition~\ref{prop:normg} is the analogue of  Entov's result that the
{\em size} of $P$ (the inverse of our $g$--norm) is larger or equal to
the $K$--area of $P$. However, in Entov's paper, the Hofer norm is taken
to be $\max_M | H |$ with $H$ normalised with $\om$--integral zero, and
it is not obvious to pass from results concerning one norm to results
concerning the other.

\end{document}